\documentclass{gtart_h}  

\def\ifplaintex{\expandafter\ifx\csname documentclass\endcsname\relax}

\def\ifplaintex{\expandafter\ifx\csname documentclass\endcsname\relax}

\ifplaintex 
\else
\fi

\expandafter\ifx\csname epsfbox\endcsname\relax\input epsf\fi

\def\gt{{\mathsurround=0pt\it $\cal G\mskip-2mu$eometry \&\ 
$\cal T\!\!$opology}}        

\def\gtp{{\mathsurround=0pt\it $\cal G\mskip-2mu$eometry \&\ 
$\cal T\!\!$opology $\cal P\!$ublications}}  

\def\lognumber#1{\def\thelognumber{#1}}
\def\volumenumber#1{\def\thevolumenumber{#1}}
\def\papernumber#1{\def\thepapernumber{#1}}
\def\volumeyear#1{\def\thevolumeyear{#1}}

\def\pagenumbers#1#2{\def\startpage{#1}\def\finishpage{#2}}
\def\published#1{\def\publishdate{#1}}
\def\proposed#1{\def\theproposer{#1}}
\def\seconded#1{\def\theseconders{#1}}
\def\received#1{\def\receiveddate{#1}}

\def\accepted#1{\def\accepteddate{#1}}

\def\asciiaddress#1{\def\theasciiaddress{#1}}
\def\asciiemail#1{\def\theasciiemail{#1}}

\long\def\asciiabstract#1{\long\def\theasciiabstract{#1}}

\def\shortauthors#1{\def\theshortauthors{#1}}

\let\\\par\let\thelognumber\relax
\let\thevolumenumber\relax\let\thepapernumber\relax
\let\thevolumeyear\relax\let\thesamplenumber\relax\let\startpage\relax
\let\finishpage\relax\let\publishdate\relax\let\receiveddate\relax
\let\reviseddate\relax\let\accepteddate\relax\let\theasciititle\relax
\let\theasciiauthors\relax\let\theasciiaddress\relax
\let\theasciiabstract\relax
\let\theasciiemail\relax\let\theshortauthors\relax\let\theshorttitle\relax

\long\def\maketitlep{   

\count0=\startpage

\gt\hfill      
\hbox to 77pt{\vbox to 0pt{\vglue -15pt\epsfbox{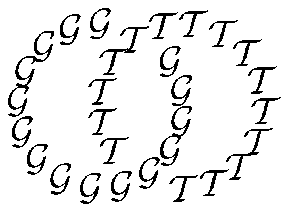}\vss}\hss}
\break
{\small\ifx\thesamplenumber\relax 
Volume \else Sample
\fi\thevolumenumber\ (\thevolumeyear)
\startpage--\finishpage\nl
Published: \publishdate}
\vglue 0.5truein plus 0.4fil minus 0.1truein

{\parskip=0pt\leftskip 0pt plus 1fil\def\\{\par\smallskip}{\ifplaintex\large
\else\Large\fi\bf\thetitle}\par\medskip}   

\vglue 0pt plus 0.1fil 

{\parskip=0pt\leftskip 0pt plus 1fil\def\\{\par}{\sc\theauthors}
\par\medskip}

\vglue 0pt plus 0.1fil 

{\small\parskip=0pt\let\newline\\
{\leftskip 0pt plus 1fil\def\\{\par}{\sl\theaddress}\par}
\expandafter\ifx\theemail\relax    
\relax\else\vglue 5pt plus 0.02fil minus 2pt\def\\{\stdspace{\rm 
and}\stdspace} 
\cl{Email:\stdspace\tt\theemail}\fi
\ifx\theurl\relax                  
\relax\else\vglue 5pt plus 0.02fil minus 2pt\def\\{\stdspace{\rm 
and}\stdspace}
\cl{URL:\stdspace\tt\theurl}\fi\par}

\vglue 7pt plus 0.3fil minus 3pt

{\bf Abstract}
\vglue 5pt plus 0.1fil minus 2pt

\theabstract

\vglue 7pt plus 0.3fil minus 3pt

{\bf AMS Classification numbers}\quad Primary:\quad \theprimaryclass

Secondary:\quad \thesecondaryclass

\vglue 5pt plus 0.3fil minus 2pt

{\bf Keywords:}\quad \thekeywords

\vglue 10pt plus 0.5fil minus 5pt

{\small  Proposed: \theproposer\hfill Received: \receiveddate\nl
Seconded: \theseconders\hfill 
\ifx\reviseddate\relax                         
Accepted: \accepteddate                        
\else
Revised: \reviseddate                          
\fi}
\eject
}       

\font\phead=cmsl9 scaled 950
\font=cmsl9 scaled 1050
\font\pnum=cmbx10 scaled 913
\font\lnum=cmbx10 
\font\pfoot=cmsl9 scaled 950
\font=cmsl9 scaled 1050
\ifplaintex
\headline{\vbox to 0pt{\vskip -4.5mm\line{\small\phead\ifnum
\count0=\startpage ISSN 1364-0380 (on line)
1465-3060 (printed) \hfill {\pnum\folio}\else\ifodd\count0\def\\{ }%
\ifx\theshorttitle\relax\thetitle\else\theshorttitle\fi\hfill{\pnum\folio}
\else\def\\{ and }{\pnum\folio}\hfill\ifx\theshortauthors\relax\theauthors
\else\theshortauthors\fi\fi\fi}\vss}}
\footline{\vbox to 0pt{\vglue 0mm\line{\small\pfoot\ifnum\count0=\startpage
\copyright\ \gtp\hfill\else
\gt, Volume \thevolumenumber\ (\thevolumeyear)\hfill\fi}\vss
}}
\else
\makeatletter
\def\@oddhead{{\small\lhead\ifnum\count0=\startpage ISSN 1364-0380 (on line)
1465-3060 (printed) \hfill {\lnum\number\count0}\else\ifodd\count0
\def\\{ }\ifx\theshorttitle\relax \thetitle \else\theshorttitle\fi\hfill
{\lnum\number\count0}\else\def\\{ and }{\lnum\number\count0}
\hfill\ifx\theshortauthors\relax 
\theauthors\else\theshortauthors\fi\fi\fi}}\def\@evenhead{\@oddhead}
\def\@oddfoot{\small\lfoot\ifnum\count0=\startpage\copyright\ \gtp\hfill\else
\gt, Volume \thevolumenumber\ (\thevolumeyear)\hfill\fi}
\def\@evenfoot{\@oddfoot}
\makeatother
\fi

\newwrite\gtoutfile
\long\gdef\makeheadfile{  
{\def\\{, }\def\s{ }
\immediate\openout\gtoutfile head.xxx
\immediate\write\gtoutfile{Proxy-for: \ifx\theasciiauthors\relax
\theauthors\else\theasciiauthors\fi\s<\ifx\theasciiemail\relax\theemail\else\theasciiemail\fi>}
\immediate\write\gtoutfile{\noexpand\\}
\immediate\write\gtoutfile{Authors: \ifx\theasciiauthors\relax
\theauthors\else\theasciiauthors\fi}
{\def\\{ }\immediate\write\gtoutfile{Title: \ifx\theasciititle\relax
\thetitle\else\theasciititle\fi}}
\immediate\write\gtoutfile{Subj-class: GT or SG or MG etc}
\immediate\write\gtoutfile{MSC-class: \theprimaryclass\ifx\thesecondaryclass\relax\else, \thesecondaryclass\fi}
\immediate\write\gtoutfile{Journal-ref: Geom. Topol. \thevolumenumber
(\thevolumeyear) \startpage-\finishpage}
\immediate\write\gtoutfile{Comments: Published by Geometry and Topology at}
\immediate\write\gtoutfile{\s\s http://www.maths.warwick.ac.uk/gt/GTVol\thevolumenumber/paper\thepapernumber.abs.html}
\immediate\write\gtoutfile{\noexpand\\}
\immediate\write\gtoutfile{}
\ifx\theasciiabstract\relax
\immediate\write\gtoutfile{\theabstract}\else
\immediate\write\gtoutfile{\theasciiabstract}\fi
\immediate\write\gtoutfile{}
\immediate\write\gtoutfile{\noexpand\\}
\immediate\write\gtoutfile{}
\immediate\closeout\gtoutfile}}  

\def\maketitlepage{\maketitlep\makeheadfile}
\let\maketitle\maketitlepage

\lognumber{428}

\volumenumber{8}\papernumber{30}\volumeyear{2004}
\pagenumbers{1079}{1125}
\received{3 March 2004}
\published{22 August 2004}
\accepted{17 July 2004}
\proposed{Haynes Miller}
\seconded{Thomas Goodwillie, Steven Ferry}

\usepackage{amscd,amssymb,amsmath,amsthm}
\usepackage[curve,knot]{xy}

\def\S{Section }

\theoremstyle{plain}
\newtheorem{theorem}[subsection]{Theorem}
\newtheorem{lemma}[subsection]{Lemma}
\newtheorem{prop}[subsection]{Proposition}
\newtheorem{cor}[subsection]{Corollary}

\newtheorem{thm}{Theorem}

\theoremstyle{definition}
\newtheorem{remark}[subsection]{Remark}
\newtheorem{definition}[subsection]{Definition}
\newtheorem{example}[subsection]{Example}
\newtheorem{question}[subsection]{Question}

\newcommand{\norm}[1]{\lVert #1 \rVert}
\renewcommand{\o}[1]{\bar{#1}}
\newcommand{\dual}[1]{{^{\sharp}#1}}
\newcommand{\comp}[1]{{#1{\,\widehat{\:\,}}}}
\newcommand{\compp}[1]{\widehat{#1}}

\renewcommand{\a}{\alpha}

\newcommand{\A}{\mathcal{A}}
\newcommand{\B}{\mathcal{B}}
\newcommand{\CC}{\mathcal{C}}
\newcommand{\HH}{\mathcal{H}}
\newcommand{\CL}{\mathcal{L}}

\newcommand{\C}{\mathbb{C}}
\newcommand{\F}{\mathbb{F}}
\newcommand{\Z}{\mathbb{Z}}
\newcommand{\Q}{\mathbb{Q}}
\newcommand{\R}{\mathbb{R}}

\newcommand{\BL}{\mathbb{L}}

\DeclareMathOperator{\rank}{rank}

\DeclareMathOperator{\id}{id}

\DeclareMathOperator{\Hom}{Hom}

\DeclareMathOperator{\ab}{ab}

\DeclareMathOperator{\im}{im}
\DeclareMathOperator{\gr}{gr}
\DeclareMathOperator{\Hilb}{Hilb}

\DeclareMathOperator{\lk}{lk}
\DeclareMathOperator{\ad}{ad}

\DeclareMathOperator{\Tor}{Tor}
\DeclareMathOperator{\coalg}{coalg}

\newcommand{\surj}{\twoheadrightarrow}

\newcommand{\cga}{\textsf{cga}}
\newcommand{\dga}{\textsf{dga}}
\newcommand{\glie}{\textsf{glie}}
\newcommand{\grlie}{\textsf{grlie}}
\newcommand{\bglie}{\textsf{bglie}}
\newcommand{\clie}{\textsf{clie}}
\newcommand{\mlie}{\textsf{mlie}}
\newcommand{\dgl}{\textsf{dgl}}
\newcommand{\cgc}{\textsf{cgc}}

\newcommand{\DGA}{\mathbf{DGA}}
\newcommand{\DGL}{\mathbf{DGL}}

\renewcommand{\H}{\operatorname{Hopf}}
\newcommand{\Der}{\operatorname{Der}_n}
\newcommand{\wL}{\widehat{\BL}}
\newcommand{\G}{\mathsf{G}}

\begin{document}

\title{Homotopy Lie algebras, lower central series\\and the Koszul property}

\authors{Stefan Papadima\\Alexander I Suciu} 
\shortauthors{S Papadima and A Suciu}

\address{Institute of Mathematics of the Romanian Academy\\PO
Box 1-764, RO-014700 Bucharest, Romania}
\secondaddress{Department of Mathematics, Northeastern 
University\\Boston, MA 02115, USA}

\asciiaddress{Institute of Mathematics of the Romanian Academy\\PO
Box 1-764, RO-014700 Bucharest, Romania\\and\\Department of 
Mathematics, Northeastern University\\Boston, MA 02115, USA}

\asciiemail{stefan.papadima@imar.ro, a.suciu@neu.edu}
\gtemail{\mailto{stefan.papadima@imar.ro}\qua {\rm and}\qua 
\mailto{a.suciu@neu.edu}}




\begin{abstract}
Let $X$ and $Y$ be finite-type CW--complexes ($X$ connected,
$Y$ simply connected), such that the rational cohomology ring of $Y$
is a $k$--rescaling of the rational cohomology ring of $X$.
Assume $H^*(X,\Q)$ is a Koszul algebra.
Then, the homotopy Lie algebra $\pi_*(\Omega Y)\otimes \Q$
equals, up to $k$--rescaling, the graded rational Lie algebra
associated to the lower central series of $\pi_1(X)$.
If $Y$ is a formal space, this equality
is actually equivalent to the Koszulness of $H^*(X,\Q)$.
If $X$ is formal (and only then), the equality lifts to  a filtered
isomorphism between the Malcev completion of $\pi_1(X)$
and the completion of $[\Omega S^{2k+1},\Omega Y]$.
Among spaces that admit naturally defined homological
rescalings are complements of complex hyperplane
arrangements, and complements of classical links. The
Rescaling Formula holds for supersolvable arrangements,
as well as for links with connected linking graph.
\end{abstract}

\asciiabstract{Let X and Y be finite-type CW--complexes (X connected,
Y simply connected), such that the rational cohomology ring of Y
is a k-rescaling of the rational cohomology ring of X.
Assume H^*(X,Q) is a Koszul algebra.
Then, the homotopy Lie algebra pi_*(Omega Y) tensor Q
equals, up to k-rescaling, the graded rational Lie algebra
associated to the lower central series of pi_1(X).
If Y is a formal space, this equality
is actually equivalent to the Koszulness of H^*(X,Q).
If X is formal (and only then), the equality lifts to  a filtered
isomorphism between the Malcev completion of pi_1(X)
and the completion of [Omega S^{2k+1} ,Omega Y].
Among spaces that admit naturally defined homological
rescalings are complements of complex hyperplane
arrangements, and complements of classical links. The
Rescaling Formula holds for supersolvable arrangements,
as well as for links with connected linking graph.}

\maketitlepage

\section{Introduction and statement of results}
\label{sec:intro}

\subsection{A tale of two Lie algebras}
\label{subsec:liealg}
In this paper, we focus on two Lie algebras, traditionally
associated to a topological space.
All spaces under consideration are connected, well-pointed, and have
the homotopy type of a finite-type CW--complex (for short, finite-type 
CW--spaces). We use rational coefficients, unless otherwise specified.

Let $X$ be a connected CW--space of finite type, with fundamental group 
$G=\pi_1(X)$. The {\em associated graded Lie algebra}\/ of $X$ is by definition
\begin{equation}
\label{eq:assocgrlie}
L_*(X):= \bigoplus_{n\geq 1} \gr^n(\pi_1(X))\otimes \Q.
\end{equation}
Here $\gr^n (G)=\Gamma_nG/\Gamma_{n+1}G$, where $\{\Gamma_n G\}_n$ denotes the
{\em lower central series} of $G$,  inductively defined by $\Gamma_1G=G$, and
$\Gamma_{n+1}G=(G,\Gamma_n G)$, the subgroup generated by all commutators
$xyx^{-1}y^{-1}$, with $x\in G$ and $y\in  \Gamma_n G$.
The graded vector space $L_*(X)=\gr^*(\pi_1(X))\otimes \Q$ has the
structure of a {\em Lie algebra with grading} (\grlie), which means
that the Lie bracket induced by the group commutator is of degree zero,
and the usual Lie identities are satisfied.

Now let $Y$ be a simply-connected CW--space of finite type.  The {\em
(rational) homotopy Lie algebra} of $Y$ is by definition
\begin{equation}
\label{eq:homotopylie}
E_*(Y):=\bigoplus_{n\geq 1}\pi_n(\Omega Y)\otimes \Q.
\end{equation}
The Whitehead product on $\pi_*(Y)$
defines, via the boundary map in the path fibration over $Y$, a
Lie bracket on $\pi_*(\Omega Y)$, called the
Samelson product (see eg~\cite{Wh}). This
makes $E_*(Y)=\pi_*(\Omega Y)\otimes \Q$
into a {\em graded Lie algebra} (\glie), meaning that the
Lie bracket is of degree zero, and the Lie identities
are satisfied only up to sign, following the Koszul convention.

\subsection{Rescaling}
\label{subsec:rescale}
We now define the basic operation of {\em rescaling},
on algebras and Lie algebras, following~\cite{CCX}.

Let $A^*$ be a connected, graded, graded-commutative, $\Q$--algebra (\cga).
For each $k\geq 0$, we define a new \cga, denoted by $A[k]$ and
called the {\em $k$--rescaling} of $A$, by
\begin{equation}
\label{eq:algresc}
A[k]^p =
\begin{cases}
A^q &\quad\text{if $p=q(2k+1)$,}\\
0   &\quad\text{otherwise,}
\end{cases}
\end{equation}
with multiplication
$A[k]^{q(2k+1)}\otimes A[k]^{q'(2k+1)}\to A[k]^{(q+q')(2k+1)}$
identified with the multiplication $A^{q}\otimes A^{q'}\to A^{q+q'}$. 
The Hilbert series of the two graded algebras are related, when defined, 
by $\Hilb(A[k],t)=\Hilb(A,t^{2k+1})$.

Let $L_*$ be a Lie algebra with grading.
For each $k\geq 1$, we define
a new \grlie, denoted by $L[k]$ and called 
the {\em $k$--rescaling} of $L$, by
\begin{equation}
\label{eq:lieresc}
L[k]_p = \begin{cases}
L_q &\quad\text{if $p=2kq$,}\\
0   &\quad\text{otherwise,}
\end{cases}
\end{equation}
with Lie bracket
$L[k]_{2kq}\otimes L[k]_{2kq'}\to L[k]_{2k(q+q')}$
identified with the Lie bracket $L_{q}\otimes L_{q'}\to L_{q+q'}$.
Note that $L[k]$ is evenly graded, and therefore may be also viewed as a
graded Lie algebra. 
When defined, $\Hilb(L[k],t)=\Hilb(L,t^{2k})$.

\begin{definition}
Let $X$ be a connected, finite-type CW--space,
and $k$ a positive integer.
A simply-connected, finite-type CW--space $Y$ is called a
({\em homological}) {\em $k$--rescaling}\/ of $X$ if
\begin{equation}
\label{eq:hxy}
H^*(Y,\Q)=H^*(X,\Q)[k].
\end{equation}
\end{definition}

Our main goal in this paper is to understand when the homological 
rescaling property \eqref{eq:hxy} passes to homotopy groups and implies
the following ({\em homotopy}) {\em Rescaling Formula}:
\begin{equation}
\label{eq:RF}
E_*(Y) = L_*(X)[k]\, .
\end{equation}

\subsection{The work of Cohen--Cohen--Xicot\'{e}ncatl}
\label{subsec:CCX}
One instance where the above formula holds is provided by
the setup considered in \cite{CCX} (see also \cite{CX}).
Let $\A=\{H_1,\dots,H_n\}$ be an arrangement of hyperplanes 
in a complex vector space $V$. For each $k\geq 1$, 
Cohen, Cohen and Xicot\'{e}ncatl construct a
so-called {\em redundant} arrangement,
$\A^k:=\{H_1^{\times k},\dots,H_n^{\times k}\}$, of codimension
$k$  subspaces in $V^{\times k}$. Their motivation was to better 
understand, in this way, previously discovered relations between 
the topology of the configuration spaces of $\C$, and the topology 
of the corresponding higher configuration spaces
(Cohen--Gitler~\cite{CG1}, Fadell--Husseini~\cite{FH}).
Indeed, if one starts with the braid arrangement, $\B_{\ell} =
\{z_i-z_j=0\mid i<j\}$, in $\C^{\ell}$,
with complement $M(\B_{\ell})$ equal to the configuration space
of $\ell$ distinct points in $\C$, then it is immediate to see that
$M(\B_{\ell}^k)$ is the configuration space of
$\ell$ distinct points in $\C^k$.

Now let $\A$ be an arbitrary hyperplane arrangement, with
complement $X=M(\A)$.  Fix an integer $k\ge 1$, and set
$Y=M(\A^{k+1})$. As shown in
\cite[Corollary~2.3]{CCX}, the homological rescaling
property \eqref{eq:hxy} is satisfied by $X$ and $Y$,
even with $\Z$ coefficients.
Moreover, Theorem 1.3 from \cite{CCX}
implies that the homotopy rescaling formula \eqref{eq:RF}
holds as well, provided $\A$ is a {\em fiber-type} arrangement.

Essentially, the proof of \cite[Theorem 1.3]{CCX}---which also gives 
information over the integers---is based upon the following geometric idea. 
Assuming that $\A$ is fiber-type amounts to saying that $M(\A)$ has an 
iterated (split) fibered structure, with trivial monodromy action on 
homology, where all homotopy fibers are wedges of circles. It turns 
out that $M(\A^{k+1})$ follows the same pattern, with the circles 
replaced by $(2k+1)$--spheres, whence the homotopy rescaling.

\subsection{Formality properties}
\label{subsec:formal}
We will use a different approach, based on rational homotopy theory, 
and D Sullivan's notion of {\em formality} (see \cite{S}).
Formal spaces are characterized by the property that their 
rational homotopy type is a formal consequence of their rational 
cohomology algebra. As such, they lend themselves to various algebraic
computations which may provide valuable homotopy information. For instance,
the main result from \cite{PY}, which holds for all formal spaces $X$, 
states:  The (Bousfield--Kan) rationalization $X_{\Q}$ is a $K(\pi,1)$ 
space if  and only if $H^*(X,\Q)$ is a Koszul algebra.

A fundamental result of Sullivan explains the relationship between 
the rational cohomology algebra of a formal space $X$ and the 
associated graded Lie algebra, $L_*(X)$. More precisely, it says 
that the {\sf grlie} $L_*(X)$ is determined by the cup-product map  
in low degrees, $H^1(X, \Q)\wedge H^1(X, \Q) \to H^2(X, \Q)$
(see \cite{PS04} for metabelian versions of this result). 

Formal computations work well for arrangement complements, because of
the following two results. First, all hyperplane arrangement complements 
are formal spaces, as follows from Brieskorn's solution \cite{Br} of 
a conjecture of Arnold. Second, all redundant subspace arrangement 
complements are formal, as well.  This is a consequence of recent 
work by Yuzvinsky \cite{Y}, sharpening an earlier result of 
De~Concini and Procesi \cite{DCP}. Indeed, if the intersection 
poset of a subspace arrangement $\B$ is a geometric lattice
(which always happens for $\B=\A^k$, as noted in \cite[Proposition~2.2]{CCX}), 
then the complement $M(\B)$ is a formal space, as noted in 
\cite[Remark 7.3(ii)]{Y}.

\subsection{Koszulness and the rescaling formula}
\label{subsec:rescaleKoszul}

Let $A^*$ be a connected, graded algebra (not necessarily
graded-commutative), over a field $\Bbbk$.
By definition, $A$ is a {\em Koszul algebra}
if $\Tor_{p,q}^A(\Bbbk,\Bbbk)=0$, for all $p\ne q$.
This is a fundamental notion
in homological algebra, going back to Priddy's analysis
of the Steenrod algebra \cite{Pr}; see
Beilinson, Ginzburg and Soergel~\cite{BGS} as a basic reference.

Our first  main result is the following.

\begin{thm}
\label{thm:RF1}
Let $X$ be a connected, finite-type CW--space, and let
$Y$ be a formal, simply-connected, finite-type CW--space.
Let $L_*(X)=\gr^*(\pi_1(X))\otimes \Q$ be the associated
graded Lie algebra of $X$, and let $E_*(Y)=\pi_*(\Omega Y)\otimes \Q$ be
the homotopy  Lie algebra of $Y$.  Assume $H^*(Y,\Q)=H^*(X,\Q)[k]$,
for some $k\geq 1$. Then the following conditions are equivalent:
\begin{enumerate}
\item \label{a1}
The graded vector spaces $E_*(Y)$ and $L_*(X)[k]$ have the 
same Hilbert series.
\item \label{a2}
The graded Lie algebras $E_*(Y)$ and $L_*(X)[k]$ are isomorphic.
\item \label{a3}
$H^*(X,\Q)$ is a Koszul algebra.
\end{enumerate}
\end{thm}
\begin{remark}
\label{rem:RFvsCCX}
As shown by Shelton and Yuzvinsky~\cite{SY}, the cohomology algebra 
$H^*(M(\A),\Q)$ is Koszul, provided $\A$ is a fiber-type arrangement.  
(The converse is an open question.)
Together with the discussion from \S\ref{subsec:CCX} 
and \S\ref{subsec:formal}, this shows that the 
implication \eqref{a3} $\Longrightarrow$ \eqref{a2} in the theorem 
above is a potential strengthening of Theorem 1.3 from
\cite{CCX}, at least for $\Q$--coefficients. At the same time, 
the implication \eqref{a2} $\Longrightarrow$ \eqref{a3} 
offers a more illuminating approach to the homotopy rescaling 
formula \eqref{eq:RF}, indicating the key role played by the 
Koszul property.
\end{remark}
\begin{remark}
\label{rem:test}
There is a well-known {\em Koszul duality formula}
relating the Hilbert series of a finite type Koszul algebra to
that of its quadratic dual:  $\Hilb(A,t) \cdot \Hilb(A^{!},-t)=1$;
see \cite{BGS}.
It is also known that this equality does not imply Koszulness;
see \cite{R}. From this point of view, the equivalence
\eqref{a1} $\Longleftrightarrow$ \eqref{a3} from Theorem~\ref{thm:RF1}
may be interpreted as a seemingly new, necessary and sufficient
test for Koszulness of \cga's; see Corollary~\ref{cor:Ktest}.
\end{remark}

\begin{remark}
\label{rem:asymp}
If $H^{>d}(X, \Q)=0$ (eg, if $X$ has dimension~$d$), 
then one knows from \cite{SYa} that $Y$ must be formal, 
as soon as $2k+1>d$; see Proposition~\ref{sss1}.
\end{remark}

\begin{remark}
\label{rem:finite type}
The finite-type hypotheses from Theorem \ref{thm:RF1} (and subsequent 
results) are imposed by the topology--algebra dictionary from rational 
homotopy theory; see Bousfield and Gugenheim \cite[Remark 11.5]{BG}. 
In this context, we should mention an interesting family of examples, 
related to braids on closed oriented surfaces of positive genus, 
analyzed by Cohen--Xicot\'{e}ncatl \cite{CX} and 
Cohen--Kohno--Xicot\'{e}ncatl \cite{CKX}.  The pairs of spaces 
coming from this family satisfy the homotopy rescaling 
formula  \eqref{eq:RF}, without being of finite type.  
It is not clear whether, for these examples, one can 
derive \eqref{eq:RF} from our Theorem \ref{thm:RF1}, 
by passing to the limit over finite-type subcomplexes.
\end{remark}

The key word in Theorem~\ref{thm:RF1} is Koszulness. Our
next result shows that this property is strong enough
to derive \eqref{eq:RF} from \eqref{eq:hxy}, even without
formality assumptions.

\begin{thm}
\label{thm:RF2}
Let $X$ and $Y$ be finite-type CW--spaces ($X$ $0$--connected, and $Y$
$1$--connected), such that $H^*(Y,\Q)=H^*(X,\Q)[k]$,
for some $k\ge 1$. If $H^*(X,\Q)$ is a Koszul algebra, then the graded
Lie algebras $E_*(Y)$ and $L_*(X)[k]$ are isomorphic.
\end{thm}

In other words, the implication \eqref{a3} $\Longrightarrow$ \eqref{a2}
from Theorem~\ref{thm:RF1} does not need the formality of  $Y$.
On the other hand, the proof of Theorem~\ref{thm:RF2} relies heavily on
the formal case.

\subsection{Coformality obstruction}
\label{subsec:coformality}
There is also a notion of {\em coformality},
due to Neisendorfer and Miller \cite{NM},
and which is the Eckmann--Hilton dual to the notion of formality.
A simply-connected space is coformal if its rational homotopy type
is determined by its homotopy Lie algebra.
In view of Theorem~\ref{thm:RF1}, the next result
offers a topological obstruction to the homotopy rescaling formula.

\begin{prop}
\label{prop:coformal}
Let $X$ and $Y$ be finite-type CW--spaces ($X$ connected, $Y$ formal and
simply-connected), such that $H^*(Y,\Q)=H^*(X,\Q)[k]$, for some $k\ge 1$.
If $H^*(X,\Q)$ is a Koszul algebra, then $Y$ is a coformal space.
\end{prop}

The converse does not hold in general, see Example~\ref{ex:nonkoszul}.

If $\A$ is an affine, generic arrangement of $n$ hyperplanes
in $\C^{n-1}$ ($n>2$), then the Rescaling Formula fails
for $X=M(\A)$, due to the non-coformality of $Y=M(\A^{k+1})$;
the absence of coformality is detected by higher-order
Whitehead products, see Example~\ref{ex:threelines}.

\subsection{LCS formula for higher homotopy groups}
\label{subsec:highLCS}

Let $X$ be a connected CW--space of finite type. Set
$\phi_n :=\rank \gr^n (\pi_1(X))$, for $n\ge 1$.
Let $P_{X}(t)=\Hilb(H^*(X;\Q),t)$ be the Poincar\'{e} series of $X$.
The following {\em lower central series (LCS) formula}
has received considerable attention:
\begin{equation}
\label{eq:lcs}
\prod_{n\ge 1} (1-t^n)^{\phi_n}=P_{X}(-t) \, .
\end{equation}
This formula was established for classifying spaces of pure braids
by Kohno \cite{Ko}, and then for complements of
arbitrary fiber-type arrangements by Falk and Randell
\cite{FR}. A variant of the LCS formula holds for the 
more general class of hypersolvable arrangements, cf~\cite{JP}.

The LCS formula was related to Koszul duality in \cite{SY} 
and \cite{PY}.  Let $A$ be a connected, finite-type graded-commutative 
algebra, with associated holonomy Lie algebra $\HH_*(A)$, 
defined in \S\ref{subsec:Chen} below. Set $\phi_n:= \dim_{\Q}\HH_n(A)$, 
for $n\ge 1$. If $A$ is a Koszul algebra, then
\begin{equation}
\label{eq:klcs}
\prod_{n\ge 1} (1-t^n)^{\phi_n}=\Hilb(A,-t)\, .
\end{equation}

Our next result gives an LCS-type formula for the rational 
higher homotopy groups.  For a $1$--connected, finite-type 
CW--space $Y$, and an integer $n\ge 1$, set
$\Phi_n:=\rank \pi_n(\Omega Y)$.

\begin{thm}
\label{thm:HLCS}
Let $Y$ be a simply-connected CW--space of finite type. Assume
that $H^*(Y, \Q)$ is the $k$--rescaling of a Koszul algebra.
Then $\Phi_n =0$, if $n$ is not
a multiple of $2k$, and the following
{\em homotopy LCS formula} holds:
\begin{equation}
\label{eq:hlcs}
\prod_{r=1}^{\infty} (1-t^{(2k+1)r})^{\Phi_{2kr}}=P_{Y}(-t)\, .
\end{equation}
\end{thm}

\subsection{Malcev completions and coalgebra maps}
\label{subsec:mcomp}

Let $Y$ be a based, simply-connected CW--space of finite type.
For an arbitrary connected, based CW--complex $K$, let $[K,\Omega Y]$ 
denote the group (under composition of loops) of pointed homotopy
classes of based maps.

Now assume $K$ is of finite type. Then
$K$ can be filtered by an increasing sequence
of connected, finite subcomplexes, $\{K_r\}_{r\ge 0}$, with
$K_0=\text{base point}$, and $\bigcup_{r\ge 0} K_r=K$.
The group $[K,\Omega Y]$ inherits a natural filtration,
with $r$-th term equal to
$\ker\, ([K, \Omega Y] \surj [K_{r-1},\Omega Y])$.
As noted by Cohen and Gitler in \cite{CG2},
the filtered group $[\Omega S^2,\Omega Y]$
is a particularly interesting object. As a set, it equals
$\prod _{n\ge 1} \pi_n(\Omega Y)$, thus reassembling all
the homotopy groups of $Y$ into a single, naturally
defined group (the ``group of homotopy groups" of $\Omega Y$).

In general, the filtered groups $[K,\Omega Y]$ are very
difficult to handle. Even so, there are two analogues
which are often easier to compute.

The first one is the Milnor--Moore group of degree~$0$ coalgebra maps,
\begin{equation}
\label{eq:mm}
\Hom^{\coalg} (H_*(K,\Q),H_*(\Omega Y,\Q)),
\end{equation}
defined in \cite[\S8]{MM}.
(When $Y$ is a $\Q$--local space, this group is in fact
isomorphic to $[K,\Omega Y]$; see \cite{Sc}.)
There is a natural filtration on
the group \eqref{eq:mm}, with $r$-th term equal to
\begin{equation}
\label{eq:filt1}
\ker\big( \Hom^{\coalg} (H_*(K),H_*(\Omega Y))\longrightarrow
\Hom^{\coalg} (H_*(K_{r-1}),H_*(\Omega Y)) \big),
\end{equation}
the kernel of the map induced by the inclusion $K_{r-1}\to K$; 
see \cite{CG2}.

The second one is the group $\comp{[K,\Omega Y]}:=
\varprojlim_r ([K_{r-1},\Omega Y]\otimes\Q)$,
which comes endowed with the canonical
inverse limit filtration, with $r$-th term equal to
\begin{equation}
\label{eq:filt2}
\ker\big( \comp{[K,\Omega Y]} \longrightarrow 
[K_{r-1},\Omega Y]\otimes \Q\big).
\end{equation}

In the above construction, we used the Malcev completion 
of a group $G$, denoted by $G\otimes \Q$; see 
Quillen~\cite[Appendix~A]{Q}.  The group
$G\otimes \Q=\varprojlim_n ( (G/\Gamma_n G)\otimes \Q )$
comes equipped with the canonical (Malcev) inverse limit filtration.

If $K=\Omega S^m$, $m\ge 2$, then $K$ has a cell decomposition
(due to M.~Morse and I. James), with one cell of dimension $r(m-1)$,
for each $r\ge 0$ (see \cite{Mi} and \cite{Wh}). 
Setting $K_r$ equal to the $r(m-1)$-th skeleton, we obtain,
by the above procedure, the filtered group
$\comp{[\Omega S^m,\Omega Y]}$.

\begin{thm}
\label{thm:MF}
Let $X$ and $Y$ be finite-type CW--spaces ($X$ connected and
$Y$ simply-connected), such that $H^*(Y,\Q)=H^*(X,\Q)[k]$,
for some $k\ge 1$. Assume that $H^*(X,\Q)$ is a Koszul algebra.
Then the next two properties are equivalent:
\begin{enumerate}

\item \label{Db}
$X$ is a formal space.

\item \label{Da}
The following filtered groups are isomorphic:

\begin{enumerate}

\rm\item \label{gr1}\sl
$\comp{[\Omega S^{2k+1}, \Omega Y]}$, with the 
filtration \eqref{eq:filt2};

\rm\item \label{gr2}\sl
$\Hom^{\coalg} (H_*(\Omega S^{2k+1},\Q),H_*(\Omega Y,\Q))$,
 with the filtration \eqref{eq:filt1};

\rm\item \label{gr3}\sl
$\pi_1(X) \otimes \Q$, with the Malcev filtration.
\end{enumerate}

\end{enumerate}
\end{thm}

This theorem lifts the Rescaling Formula \eqref{eq:RF} from the 
level of associated graded Lie algebras to the level of filtered 
groups; see Remark~\ref{rem:upgrade}.  The implication 
\eqref{Db} $\Longrightarrow$ \eqref{Da} from Theorem~\ref{thm:MF}
follows from Theorem~\ref{thm:RF2}, via a well-known
formula (valid for any formal space $X$), relating the group
$\pi_1(X)\otimes \Q$ to the Lie algebra $L_*(X)$, see \eqref{eq:f3},
together with two general results (valid for
an arbitrary\break $1$--connected, finite-type CW--space $Y$),
proved in \S\ref{sec:thmD}.

The first result (Theorem~\ref{thm:exp}) gives a description
of the filtered group of homotopy classes,
$\comp{[\Omega S^m,\Omega Y]}$, $m\ge 2$,
in terms of the homotopy Lie algebra, $E_*(Y)$.
The proof depends heavily on a formula of Baues \cite{Ba},
reviewed in \S\ref{subsec:baues}.
The Rescaling Formula may then be invoked to
establish the isomorphism between the filtered groups
\eqref{gr1} and \eqref{gr3} from Theorem~\ref{thm:MF}.

The second result (Proposition~\ref{prop:homco}) gives
a filtered isomorphism between the groups $\comp{[K,\Omega Y]}$ and
$ \Hom^{\coalg} (H_*(K),H_*(\Omega Y))$, for any connected, finite-type
CW--complex $K$, filtered as in the beginning of \S\ref{subsec:mcomp}.

\subsection{Connection with the work of Cohen--Gitler--Sato}
\label{subsec:CGS}

The two results mentioned above may be combined to obtain information
about the filtered Milnor--Moore groups
\begin{equation}
\label{eq:mms2}
\Hom^{\coalg} (H_*(\Omega S^2,\Q),H_*(\Omega Y,\Q)),
\end{equation}
studied by Cohen--Gitler \cite{CG2}, Sato~\cite{Sa}, 
and Cohen--Sato \cite{CS}.

For instance, the Lie algebra with grading associated to the
filtered group \eqref{eq:mms2} is isomorphic to a certain
rebracketing of the homotopy Lie algebra $E_*(Y)$. When $Y$ is a
$k$--rescaling of a space $X$ with Koszul cohomology algebra, this
rebracketing equals $L_*(X)[k]$. From these considerations,
we recover several results from \cite{CG2} and \cite{CS};
see Corollary~\ref{grglor} and Remark~\ref{rem:CGS2} for details.

Furthermore, we show in Proposition~\ref{prop:answer} that the 
Milnor--Moore group \eqref{eq:mms2} is isomorphic to the Malcev 
completion of  $\pi_1(X)$, provided $Y$ is a $k$--rescaling of a formal 
space $X$, satisfying the assumptions of Theorem~\ref{thm:MF}.
In the case when $X=M(\B_{\ell})$ and $Y=M(\B_{\ell}^{k+1})$,
this answers a question raised in \cite{CG2}; see Remark ~\ref{rem:CGS1}.
In the case when $X=\bigvee^n S^{1}$ and $Y=\bigvee^n S^{2k+1}$, 
this recovers a result of \cite{Sa}; see Example ~\ref{ex:malsimp}.

\subsection{Rescaling links in $S^3$}
\label{subsec:links}
Using Sullivan's minimal models,
it is easy to see that any connected CW--space of finite type, $X$,
admits a homological $k$--rescaling, for each $k\ge 1$. In general,
though, such a rescaling, $Y$, is not uniquely determined  by condition
\eqref{eq:hxy}, not even up to rational homotopy equivalence.
See \S\ref{subsec:existence} for details.

Besides complements of hyperplane arrangements, there is
another large class of spaces that admit a naturally
defined rescaling:  complements of classical links.

Let $K=(K_1,\dots,K_n)$ be a link of oriented circles in $S^3$.
Associated to $K$ there is a {\em linking graph}, $\G_K$,
with vertices corresponding to the components $K_i$,
and (simple) edges connecting pairs of distinct vertices for which
$\lk(K_i,K_j)\ne 0$.

We define $K^{\circledast k}$ (the $k$--rescaling of $K$)
to be the link of $(2k+1)$--spheres in $S^{4k+3}$ obtained by taking the
iterated join (in the sense of Koschorke and Rolfsen \cite{KR})
of the link $K$ with $k$ copies of the $n$--component Hopf link.
It turns out that the complement of $K^{\circledast k}$
is indeed the (unique up to rational homotopy)
$k$--rescaling of the complement of $K$; see Proposition~\ref{prop:reslice}.

\begin{thm}
\label{thm:RFlk}
Let $K=(K_1,\dots,K_n)$ be a link in $S^3$, with complement $X$.
Let $K^{\circledast k}$ be the $k$--rescaling of $K$, 
with complement $Y$.  Then the Rescaling Formula 
holds for $X$ and $Y$ if and only if the linking graph
of $K$ is connected:
\[
E_*(Y)\cong L_*(X)[k] \Longleftrightarrow \pi_0(\G_K)=0.
\]
\end{thm}

This theorem (together with the next corollary), will be proved
in \S\ref{subsec:pfE}.
Key to the proof is the fact that $\G_K$ is connected if and
only if $H^*(X,\Q)$ is Koszul, cf~\cite{MP}.

\begin{cor}
\label{cor:lkresc}
Let $K$ be a link in $S^3$ with connected linking graph, and
let $Y$ be the complement of $K^{\circledast k}$. Then $Y$ is
both formal and coformal, and its homotopy Lie algebra is a
semidirect product of free Lie algebras  generated in degree $2k$,
\begin{equation*}
\label{eq:semi}
E_*(Y)=\BL(x_1,\dots ,x_{n-1}) \rtimes\BL(x_n),
\end{equation*}
with ranks $\Phi_q$ vanishing if $2k \nmid q$, and
\[
\prod_{r\ge 1} (1-t^{(2k+1)r})^{\Phi_{2kr}}=(1-t^{2k+1})(1-(n-1)t^{2k+1}).
\]
Moreover, $P_{\Omega Y}(t)=\frac{1}{(1-t^{2k})(1-(n-1)t^{2k})}$.
\end{cor}

There is a rich supply of links with connected linking graph. This may be
seen by combining the following two well-known facts. First, any graph
on $n$ vertices, weighted by integers, may be realized as the weighted 
linking graph of some closed-up pure braid on $n$ strands. Second, 
two links obtained by closing up pure braids which differ by a pure 
braid commutator have the same weighted linking graph.
Geometrically defined examples of links with complete (hence, connected)
linking graphs include algebraic links and singularity links of central
arrangements of transverse planes in $\R^4$.

\subsection{Formality obstructions for classical links}
\label{subsec:final}
Let $Y$ be a $k$--rescaling of a connected space, $X$,
having the homotopy type of a finite CW--complex.
If $k$ is large enough, we know from Theorem~\ref{thm:RF1}
(and Remark~\ref{rem:asymp})  that the only
obstruction to the Rescaling Formula \eqref{eq:RF}
is the Koszul property of $H^*(X, \Q)$. If this property holds,
we know from Theorem~\ref{thm:MF} that the only obstruction
to upgrading the Rescaling Formula \eqref{eq:RF} to the
{\em Malcev Formula},
\begin{equation}
\label{eq:Mfinal}
\pi_1(X)\otimes \Q \cong \comp{[\Omega S^{2k+1}, \Omega Y]}\, ,
\end{equation}
is the formality of $X$.

When $X$ is the complement of a classical link with
connected linking graph, we present, in Corollary~\ref{cor:obstr},
a sequence of obstructions to the formality of $X$, based on the
{\em Campbell--Hausdorff invariants} of links, introduced in \cite{P97}.
This leads to examples which show that, in general, the Koszulness of
$H^*(X,\Q)$ alone does not suffice to imply the Malcev Formula
\eqref{eq:Mfinal}; see Example~\ref{ex:twisthopf}.

\subsection{Organization of the paper}
\label{subsec:org}

The paper is divided into three parts, each one subdivided into
three sections:

Part~\ref{part:RF} deals with the Rescaling Formula,
and proves Theorems~\ref{thm:RF1}, \ref{thm:RF2}, and
\ref{thm:HLCS}.

Part~\ref{part:MF} deals with the Malcev Formula and Milnor--Moore groups,
and proves Theorem~\ref{thm:MF}.

Part~\ref{part:arrlk} contains applications to arrangements and links,
and proves Theorem~\ref{thm:RFlk}.  

A more detailed guide to the contents can be found at the 
beginning of each part.  An announcement of the results 
of this paper appeared in \cite{PS}. 

\subsection{Acknowledgments}
\label{subsec:ack}
S Papadima was partially supported by CERES Grant 152/2001 
of the Romanian Ministry  of Education and Research.  
A Suciu was partially supported by NSF grant DMS-0105342.

\part{\Large\bf The rescaling formula}
\label{part:RF}

In this first part, we study homological rescalings of spaces,
and prove Theorems \ref{thm:RF1}--\ref{thm:HLCS} from the Introduction,
using the algebraic models of Quillen and Sullivan as a key tool.

The main goal is to establish the homotopy Rescaling Formula,
$E_*(Y)=L_*(X)[k]$, which allows
one to reconstruct the homotopy Lie algebra of a simply-connected
space $Y$ from the Lie algebra associated to the lower central series
of the fundamental group of a space $X$, provided  $H^*(X,\Q)$ is
a Koszul algebra, and $H^*(Y,\Q)$ is a $k$--rescaling of that algebra.

The Rescaling Formula is strong enough to imply the Koszulness 
of $H^*(X,\Q)$, in the case when $Y$ is formal.
Furthermore, the formula can be used to determine 
the rational homotopy type of $\Omega Y$, solely 
from the Poincar\'{e} polynomial of $X$.

\section{Algebraic models of spaces}
\label{sec:models}

We start by recalling some basic facts which allow, in the case
of formal spaces, to describe the Lie algebras appearing in the
Rescaling Formula in terms of cohomology algebras.

\subsection{Quillen's functors and the homotopy Lie algebra}
\label{subsec:Quillen}
Let $\DGA$ be the category of differential graded algebras (\dga's).
The objects, $(A,d)$, are finite-type \cga's, $A^*$, endowed with an
algebra differential, $d$, of degree $+1$; the morphisms are the 
$\Q$--linear maps which preserve all the existing structure. Similarly, 
let $\DGL$ be the category of differential graded Lie algebras (\dgl's).
It has objects of the form $(E, \partial)$, where $E_*$ is a finite-type 
\glie~and $\partial$ is a degree $-1$ Lie differential; the morphisms 
are again required to preserve all the structure. For both categories, 
{\em quasi-isomorphisms} are morphisms inducing isomorphisms in homology, 
and {\em weak equivalences} are finite compositions (in both directions) 
of quasi-isomorphisms.

A fundamental tool in the rational homotopy theory of simply-connected spaces is
provided by Quillen's pair of adjoint functors, which relate differential
graded Lie algebras and coalgebras; see \cite{Q}. For our purposes here, we will
need a dual version; namely, the adjoint functors,
\[
\CC \co \DGL_1 \longrightarrow \DGA_1 \quad {\rm and} \quad
\CL \co \DGA_1\longrightarrow
\DGL_1\, ,
\]
between the full subcategories consisting of simply-connected objects.\ \ 
(A \dga \ $(A,d)$ is $1$--connected if $A^1=0$; a \dgl~$(E, \partial)$ is
$1$--connected if $E_*$ is strictly positively graded.)
We  briefly recall the construction of Quillen's functors,
and some of their relevant properties
\footnote{We will abuse notation, and write $\CC(E,\partial)$ for
$\CC((E,\partial))$, etc.}.
For full  details, we refer to Tanr\'{e}~\cite[Ch.~I]{Ta}.

First, we need some notation.  Let $V=V^*$ be a graded vector space over 
$\Q$. For each integer $r$, denote by $s^r V$ the $r$-th suspension of $V$,
with grading $(s^r V)^k=V^{k-r}$. Let $^{\sharp}V=\Hom_{\Q}(V,\Q)$ be the
dual of $V$, with grading $(^{\sharp}V)_k=\Hom_{\Q}(V^k,\Q)$.
Finally, let $\bigwedge V$ be the free \cga~ generated by $V$,
and let $\BL(V)$ be the free \glie~generated by $V$.
Note that $\bigwedge V$  has an
additional grading, given by word length, and denoted by
$\bigwedge^*(V) := \bigoplus_{q \ge 0} \bigwedge^q (V)$.
Similarly, $\BL(V)$
has an extra grading, given by bracket length, and denoted by
$\BL^*(V) := \bigoplus_{q \ge 1} \BL^q(V)$.

The functor $\CC \co \DGL_1 \to \DGA_1$ is defined as follows.
Let $(E,\partial)$ be a simply-connected \dgl. Then:
\begin{equation}
\label{eq:cfunctor}
\CC (E,\partial)=(\bigwedge(s\dual{E}),d),
\end{equation}
where the restriction of $d$ to the free algebra generators is of the form
$d=d_1+d_2$. By duality, the linear part,
$d_1\co s\dual{E} \to s\dual{E}$, corresponds to the Lie differential,
$\partial\co E\to E$. Likewise, the quadratic part, $d_2 \co s\dual{E}
\to s\dual{E} \wedge s\dual{E}$, comes from the Lie bracket,
$[\, , ]\co E\wedge E\to E$. Furthermore,
a \dgl~map $\varphi\co (E,\partial) \to (E',\partial')$ 
is sent to the \dga~ map
$\CC (\varphi)\co (\bigwedge s\dual{E}', d')\to (\bigwedge s\dual{E},
d)$, defined as the free \cga~ extension of the dual of $\varphi$. 
One knows that $\varphi$ is a quasi-isomorphism if and only if 
$\CC(\varphi)$ is a quasi-isomorphism.

The functor $\CL \co \DGA_1 \to \DGL_1$ is defined as follows.
Let $(A,d)$ be a simply-connected \dga. Denote by $\o{A}$ 
its augmentation ideal.  Then:
\begin{equation}
\label{eq:lfunctor}
\CL (A, d) = \big(\BL(s^{-1} \dual{\o{A}}), \partial\big),
\end{equation}
where the restriction of $\partial$ to the free Lie generators 
is of the form $\partial =\partial_1 +\partial_2$. As before, 
the linear part $\partial_1$ corresponds by duality to $d$, and the (Lie)
quadratic part $\partial_2$ comes from the dual of the algebra multiplication.
A \dga~map, $\psi\co (A,d){\to} (A',d')$, is sent to the \dgl~ map
$\CL(\psi)\co (\BL(s^{-1} \dual{\o{A}'}), \partial')
{\to} (\BL(s^{-1} \dual{\o{A}}),$ $\partial)$, defined as the free 
\glie~extension of the dual of $\psi$. Again, $\psi$ is a 
quasi-isomorphism if and only if $\CL(\psi)$ is a quasi-isomorphism.

The functors $\CC$ and $\CL$ are weakly adjoint, ie,
there exist adjunction morphisms,
\[
\alpha\co \CC\CL(A,d)\longrightarrow (A,d) \quad {\rm and} \quad
\beta\co \CL\CC(E,\partial)\longrightarrow (E,\partial)\, ,
\]
inducing isomorphisms in homology.
From the definitions, the underlying \cga~of
$\CC\CL(A,d)$ is $\bigwedge(s ^{\sharp}\BL^* (s^{-1} \dual{\o{A}}))$.
The morphism $\alpha$ is given on generators by
\begin{equation}
\label{eq:adj}
\left. \alpha\right| _{s \dual{\BL^{>1}}} = 0 \, ,\quad
\left. \alpha\right| _{s \dual{\BL^{1}}} = \id_{\o A} .
\end{equation}
There is a similar description for $\beta$, which we won't use. See
\cite[Ch.~I]{Ta} for details.

Coming back to our setup from Theorem~\ref{thm:RF1},
let us record the following first key observation.
 Let $Y$ be $1$--connected and of finite type, with
cohomology algebra $B^*=H^*(Y,\Q)$. If $Y$ is formal,
then it follows from \cite{Q} that
\begin{equation}
\label{eq:first}
E_*(Y) = H_*\CL(B,0)\, ,
\end{equation}
as graded Lie algebras.

\subsection{Chen's holonomy Lie algebra and the associated graded Lie algebra}
\label{subsec:Chen}
Let $A^*$ be a (connected) \cga~of finite type. Set $A_q :=\dual{A^q}$,
for $q \ge 1$. Define the reduced diagonal, $\nabla\co A_2 \to A_1\wedge A_1$,
to be the dual of the algebra multiplication, $\mu\co A^1\wedge A^1 \to A^2$.
The {\em holonomy (graded) Lie algebra} of $A$ is defined by:
\begin{equation}
\label{eq:hol}
\HH_*(A):= \BL^*(A_1)/\, {\rm ideal}\, (\im \nabla)\, ,
\end{equation}
with grading induced by bracket length. (We have used the standard
identification between $\BL^2(W)$ and $W\wedge W$.)

Here is our second key observation. Let $X$ be a formal, connected space
of finite type,  with cohomology algebra
$A^*=H^*(X,\Q)$. It is well-known (see for instance \cite{MP}) that
\begin{equation}
\label{eq:second}
L_*(X) = \HH_*(A)\, ,
\end{equation}
as Lie algebras with grading.  For more on the holonomy 
Lie algebra, and on its derived series quotients, see \cite{PS04}. 

\subsection{The bigraded $1$--model and the Koszul property}
\label{subsec:Koszul}
Another ingredient for the proof of Theorem~\ref{thm:RF1} is supplied by
a connection between the Koszulness of $A$,
and a certain property of $\HH(A)$.  We review this connection,
following \cite{MP} and \cite{PY}.

Set $\HH=\HH(A)$.
The free \cga~$\bigwedge^*(\dual{\HH})$ is simply the exterior
algebra on the dual of $\HH$,
when assigning to $\dual{\HH}$ the (upper) degree $1$.
It becomes a \dga~ with additional
grading, when endowed with the classical 
differential, $d_0$, that is used in defining 
Lie algebra cohomology:
\[
d_0\omega(x_1\wedge \cdots \wedge x_q)=\sum_{1\le i<j\le q} (-1)^{i+j}
\omega([x_i,x_j]\wedge x_1\wedge \cdots \wedge \hat{x}_i
\wedge \cdots \wedge \hat{x}_j \wedge \cdots \wedge x_q).
\]
Note that
$(\bigwedge(\dual{\HH}), d_0)$ is just the ungraded version of $\CC(\HH, 0)$. The
(multiplicative) extra grading is obtained by assigning to $\dual{\HH}_q$ the
(lower) degree $q-1$, for  $q \ge 1$.

Moreover, there is a canonical \dga~map (with respect to upper degrees),
\begin{equation}
\label{eq:rho}
\rho \co \big(\bigwedge\nolimits ^*(\dual{\HH}), d_0\big)
\longrightarrow (A^*,0)\, ,
\end{equation}
which sends $\dual{\HH}_q$ to zero, for $q>1$, and coincides 
with the canonical identification, $\dual{\HH}_1=A^1$, for $q=1$.  
This is the so-called bigraded $1$--modelling map of $A$; 
see \cite[Lemma 1.8(i)]{MP}. By \cite[Proposition~4.4]{PY}, we know 
that $\rho$ is a quasi-isomorphism if and only if $A$ is Koszul.

\subsection{Bigraded Lie algebras}
\label{subsec:bglie}

A {\em bigraded Lie algebra} (\bglie) is a graded Lie algebra with an
extra (upper) grading, which is preserved by the Lie bracket.
Here are two basic examples that we will use.

Let $E_*$ be a \glie, with lower central series $\{\Gamma_n E\}_{n}$,
inductively constructed by setting $\Gamma_1 E=E$, and
$\Gamma_{n+1} E= [E, \Gamma_n E]$.
The {\em associated bigraded Lie algebra},
$\gr^*E_*$, is defined as:
\begin{equation}
\label{eq:bglie}
\gr^*E_* := \bigoplus_{n\ge 1} \gr^{n}E_* \, ,\quad {\rm where} \quad
\gr^{n}E_* := \Gamma_n E_*/\Gamma_{n+1} E_* \, .
\end{equation}
Let $B^*$ be a $1$--connected \cga, and $\CL=\CL(B,0)=
(\BL(s^{-1}\dual{\o{B}}),\partial_2)$ the corresponding Quillen
minimal model; see \S\ref{subsec:Quillen}.
The Lie algebra $\CL$ is actually bigraded: the
second (upper) grading is given by bracket length in $\BL^*$.
Furthermore, the differential $\partial_2$ preserves both
gradings.  Hence, the homology of $\CL$ is also a \bglie,
to be denoted by $H^*_*(\CL)$.

\section{Homology, homotopy, and the Koszul property}
\label{sec:proofs}

In this section, we prove Theorems~\ref{thm:RF1}, \ref{thm:RF2}, 
and \ref{thm:HLCS}, as well as Proposition~\ref{prop:coformal} 
from the Introduction.  We start with a proof of Theorem~\ref{thm:RF1} 
in the particular case when the space $X$ is formal; the general 
case will be handled once Theorem~\ref{thm:RF2} is proved.

\subsection{Proof of Theorem~\ref{thm:RF1}, when $X$ is formal}
\label{subsec:pfRF1}
We start from the assumption that $Y$ is a $k$--rescaling
of $X$, for some fixed $k\ge 1$.
Set $A^*:=H^*(X,\Q)$ and $B^*:=H^*(Y,\Q)=A[k]$.
We know, from \eqref{eq:first} and \eqref{eq:second}, that
$E_*(Y)=H_*\CL(B,0)$, and $L_*(X)[k]=\HH_*(A)[k]$. We will take a first
step and prove the equivalence of \eqref{a1} and \eqref{a2} in our theorem,
by constructing a \dgl~map, $\lambda\co \CL(B,0)\to (\HH(A)[k],0)$.

To this end, we start by noting that the underlying \glie~of $\CL(B,0)$
is freely generated by the union of all $A_q:=\dual{ A^q}$ (taken in degree
$(2k+1)q-1$), for $q\ge 1$. This follows from the construction of $\CL$
explained in \S \ref{subsec:Quillen}, given our assumption that $B=A[k]$, 
see \eqref{eq:algresc}.

At the same time, \eqref{eq:hol} and \eqref{eq:lieresc} readily imply that
the \glie~ $\HH(A)[k]$ is generated by $A_1$ (taken in degree $2k$), with
defining relations $\nabla (A_2)=0$.

Sending $A_1$ identically to $A_1$, and $A_{>1}$ to zero, 
defines a \glie~map
\begin{equation}
\label{varphi}
\lambda\co \CL(B,0)\longrightarrow (\HH(A)[k],0).
\end{equation}

\begin{lemma}
\label{lem:step1}
The map $\lambda$ commutes with differentials and 
induces a surjection in homology.
\end{lemma}
\begin{proof}
To check the first assertion, it is enough to verify
it on generators. In other words, we must show that
$\lambda\partial (A_q)=0$, for all $q\ge 1$, where
$\partial$ is the differential of $\CL$ (that is, the
dual of the multiplication of $\o A$). If
$q>2$, this follows at once, for degree reasons,
given the definition of $\lambda$. If $q=2$, we
just have to note that $\partial$ and $\nabla$
coincide on $A_2$, and $\lambda$ is the identity on
$A_1$. Finally, $\partial(A_1)=0$, by construction.
This last remark also takes care of the second
assertion of the lemma.
\end{proof}

We infer from the preceding lemma, via a dimension argument,
that both property \eqref{a1} and property
\eqref{a2} from Theorem~\ref{thm:RF1} are equivalent to the
fact that $\lambda$ is a quasi-isomorphism. We
may thus finish the proof of our theorem
(in the particular case when $X$ is formal),
by verifying the following assertion.

\begin{lemma}
\label{lem:step2}
The \dgl~map $\lambda\co \CL(B,0)\to (\HH(A)[k],0)$
is a quasi-isomorphism if and
only if $A$ is a Koszul algebra.
\end{lemma}
\begin{proof} From the general theory summarized
in \S\ref{subsec:Quillen}, we know that
$\lambda$ is a quasi-isomorphism if and only if
$\alpha \circ \CC(\lambda)\co \CC(\HH[k],0) \to (B,0)$
is a quasi-isomorphism, where $\HH:=\HH(A)$, and
$\alpha\co \CC\CL(B,0) \to (B,0)$ is the adjunction map.

The underlying \cga~ of $\CC(\HH[k],0)$ is the exterior
algebra generated by the union of all $\dual{\HH}_q$
(taken in degree $2kq+1$), for $q\ge 1$; see
\S\ref{subsec:Quillen} and \eqref{eq:lieresc}. From the discussion
in \S\ref{subsec:Koszul}, we see that 
$\CC(\HH[k],0)$ equals $(\bigwedge(\dual{\HH}), d_0)$,
modulo some degree reindexing of the generators.
At the same time, $(B,0)$ equals $(A,0)$,
modulo rescaling, by hypothesis.

It remains to identify $\alpha \circ \CC(\lambda)$ with the map
$\rho$ from \eqref{eq:rho},
in order to be able to use \cite[Proposition~4.4]{PY}, and thus finish
the proof of the lemma. This we do, by looking at the action of
$\alpha \circ\CC(\lambda)$ on the algebra generators,
$\dual{\HH}_q$, for $q \ge 1$.

Since $\lambda \co \BL(s^{-1}\dual{\o{B}}) \to \HH[k]$ 
is a Lie algebra map, it preserves  bracket degree, 
and therefore $\CC(\lambda)$ sends $\dual{\HH}_q$ to
$s\dual{ \BL}^q (s^{-1}\dual{ \o{B}})$, for all $q\ge 1$.
From the way the maps $\alpha$ and $\lambda$ were defined
in \eqref{eq:adj} and  \eqref{varphi}, we find
$\alpha \circ\CC(\lambda)(\dual{\HH}_q) =0$,
if $q>1$, and $\alpha \circ \CC(\lambda)|_{\dual{\HH}_1}=\id $.
Thus, $\alpha \circ \CC(\lambda)=\rho$.
\end{proof}

The proof of Theorem~\ref{thm:RF1}
(with formality assumptions on both $Y$ and $X$)
is now complete. By examining the proof,
we derive the following Koszulness test.

\begin{cor}
\label{cor:Ktest}
Let $A$ be a finite-type, connected \cga, and $k$ a positive integer.
Then: $A$ is Koszul
if and only if $\Hilb (H_*\CL(A[k],0),t)=\Hilb (\HH_*(A)[k],t)$.
\end{cor}

\subsection{Proof of Proposition~\ref{prop:coformal}}
\label{subsec:pfcof}
We keep the notation from \S\ref{subsec:pfRF1}. By \cite{Q},
the rational homotopy type of the formal space $Y$ is given
by the weak equivalence type of its minimal Quillen model, $\CL(B,0)$.
If $A$ is a Koszul algebra,
Lemma \ref{lem:step2} may be used to infer that
$\CL(B,0)$ has the same weak equivalence type as $(E_*(Y),0)$,
whence the coformality of $Y$.

\subsection{Proof of Theorem~\ref{thm:RF2}}
\label{subsec:pfRF2}
As before, set $A^*=H^*(X,\Q)$, $B^*=H^*(Y,\Q)$;
also,  $\HH=\HH(A)$, $\CL=\CL(B,0)$. Our
assumptions, namely $B=A[k]$ and the Koszulness 
of $A$, imply that there is a \glie~ isomorphism,
\begin{equation}
\label{eq:hlambda}
H_*(\lambda)\co H_* \CL \xrightarrow{\:\simeq\:}
\HH[k]_*
\end{equation}
(see Lemmas~\ref{lem:step1} and \ref{lem:step2}).

We will prove Theorem~\ref{thm:RF2} by showing that
$L_*(X)=\HH_*$ (as \grlie's) and $E_*(Y)=H_*\CL$ (as \glie's).
To this end, we will use two results from \cite{MP}:
Theorems A'(i) and A(i), respectively.

\begin{lemma}
\label{lem:r1}
If $A^*=H^*(X,\Q)$ is a Koszul algebra, then $L_*(X)=\HH_*(A)$ (as \grlie's).
\end{lemma}

\begin{proof}
We know from \S\ref{subsec:Koszul} that the map
$\rho \co \big(\bigwedge\nolimits ^*(\dual{\HH}), d_0\big) 
\longrightarrow (A^*,0)$ is a quasi-isomorphism.  In particular, 
the bigraded minimal model of $A$ is generated in (upper) degree $1$.  
Hence, Lemma 1.9(ii) and Theorem A'(i) from \cite{MP} give the 
desired isomorphism.
\end{proof}

To prove that $E_*(Y)=H_*(\CL)$, we  need a rigidity result, at
the associated graded level. For that, we first need the following.

\begin{lemma}
\label{lem:r2}
The bigraded Lie algebra $H^*_*(\CL)$ is generated  by $H^1_*(\CL)=A_1$.
\end{lemma}

\begin{proof}
Recall that the map $\lambda\co \CL\to \HH[k]$ restricts
to the identity on $A_1$ (see \S\ref{subsec:pfRF1}).
From the definition of the holonomy Lie algebra of $A$
(see \eqref{eq:hol}), $\HH[k]$ is generated by $A_1$.
Since $H_*(\lambda)$ is a Lie isomorphism, we are done.
\end{proof}

\begin{lemma}
\label{lem:r3}
The \bglie's $\gr^* E_*(Y)$ and $H^*_*(\CL)$ are isomorphic.
\end{lemma}

\begin{proof}
The desired isomorphism is given by Theorem A(i) from \cite{MP},
provided two conditions are satisfied.  In view of 
\cite[Lemma 1.9(i)]{MP}, we have to verify that $H^*(\CL)$ 
is generated by $H^1(\CL)$, and $B^*$ is intrinsically spherically 
generated.  The first condition follows from Lemma~\ref{lem:r2}.

To verify the second
condition, note that the algebra $A^*$ is generated by
$A^1$, as a  consequence of the Koszul property (see eg~\cite{BGS}).
Therefore, the algebra $B^*=A[k]^*$ is homogeneously generated by the 
component $B^{2k+1}=A^1$. The conclusion follows from \cite[Remark 4.8]{P}.
\end{proof}

To conclude the proof of Theorem~\ref{thm:RF2}, we need to lift
the isomorphism from Lemma~\ref{lem:r3} from the associated
graded level to the Lie algebra level.

Set $\BL=\BL(s^{-1}\dual{\o{B}})$.   It is well-known
that $E_*(Y)=H_*(\BL,\partial)$, where the \dgl~$(\BL,\partial)$
is the Quillen minimal model of $Y$.  Furthermore, the restriction
of $\partial$ to $s^{-1}\dual{\o{B}}$ is of the form $\partial=\partial_2+p$,
where $\partial_2$ is the quadratic differential of $\CL=\CL(B,0)$,
and the perturbation $p$ has the property that
\begin{equation}
\label{eq:pert}
p(s^{-1}\dual{\o{B}})\subseteq \BL^{\ge 3}
\end{equation}
(see \cite{Q}, \cite{Ta}).  Recall also from \S\ref{subsec:pfRF1}
that $\BL$ is freely generated by the union of all $A_q$,
taken in degree $(2k+1)q-1$, for $q\ge 1$.

\begin{lemma}
\label{lem:r4}
$\partial \vert_{A_q} = \partial_2 \vert_{A_q} $, for $q\le 2$.
\end{lemma}

\begin{proof}
We need to show that $p(A_q)=0$, for $q\le 2$.
Since $p$ lowers degree by $1$, we have $\deg p(A_q)\le 4k$,
if $q\le 2$. On the other hand, $\deg p(A_q)\ge 6k$, for all $q\ge 1$;
see \eqref{eq:pert}.
\end{proof}

From the Lemma (with $q=1$), we see that the identity of $A_1$
extends to a \glie~map, $\eta \co \BL(A_1) \to
H_*(\BL,\partial)$.  Again from the Lemma (with $q=2$), we
infer that $\eta$ factors through a \glie~map,
\begin{equation}
\label{eq:etasurj}
\eta \co \BL(A_1)/{\rm{ideal}} (\partial_2( A_2)) \longrightarrow
H_*(\BL,\partial)=E_*(Y).
\end{equation}

\begin{lemma}
\label{lem:r5}
The \glie~map $\eta$ from \eqref{eq:etasurj} is surjective.
\end{lemma}

\begin{proof}
It is enough to show that the Lie algebra
$H_*(\BL,\partial)$ is generated by $A_1$.  Let
\begin{equation*}
\label{eq:iota}
\iota \co A_1 \longrightarrow
H_*(\BL,\partial)/[H_*(\BL,\partial),H_*(\BL,\partial)]
\end{equation*}
be the composite of the inclusion $A_1 \to H_*(\BL,\partial)$ with
the canonical projection. Note that $\iota$
is injective.  This follows from simple bracket degree reasons,
given the fact that $\partial (\BL)\subseteq \BL^{\ge 2}$.

By Lemma~\ref{lem:r3},
the vector spaces  $\gr^1 H_*(\BL,\partial)
=H_*(\BL,\partial)/[H_*(\BL,\partial),H_*(\BL,\partial)]$
(see definition \eqref{eq:bglie}) and $H^1(\CL)=A_1$
(see Lemma~\ref{lem:r2}) are isomorphic. Consequently, the map $\iota$
identifies $A_1$ with $\gr^1 H_*(\BL,\partial)$, the vector 
space of Lie algebra generators of $H_*(\BL,\partial)$.
\end{proof}

The next Corollary ends the proof of Theorem~\ref{thm:RF2}.

\begin{cor}
\label{cor:r6}
The \glie's $E_*(Y)$ and $H_*(\CL)$ are isomorphic.
\end{cor}

\begin{proof}
By Lemma~\ref{lem:r3}, the graded Lie algebras $E_*(Y)=H_*(\BL,\partial)$
and $H_*(\CL)=\HH[k]_*=\BL(A_1)/ (\partial_2( A_2))$ are isomorphic as
graded vector spaces.  Therefore, the map $\eta$ from \eqref{eq:etasurj}
is a \glie~isomorphism.
\end{proof}

\subsection{End of proof of Theorem \ref{thm:RF1}}
\label{subsec:endA}
We may now remove the formality assumption on $X$, made in
\S \ref{subsec:pfRF1}, as follows. Given Theorem~\ref{thm:RF2} and
Lemma~\ref{lem:step2}, all we have to show is that $H_*(\lambda)$ is
an isomorphism, provided that the graded vector spaces $H_*\CL (B,0)$ and
$L_*(X)[k]$ are isomorphic. Recall that one has a \grlie~surjection,
$\HH_*(A)\surj L_*(X)$; see \cite[Proposition~3.3]{MP}. Lemma~\ref{lem:step1}
and a dimension argument combine to show that $H_*(\lambda)$ is an 
isomorphism, as needed. The proof of Theorem~\ref{thm:RF1} is now complete.

\subsection{Proof of Theorem \ref{thm:HLCS}}
\label{subsec:pfHLCS}
By assumption, $H^*(Y,\Q)=A[k]^*$, with $A^*$ Koszul.
Set $\HH_*=\HH_*(A)$.  We may write the homotopy 
Rescaling Formula \eqref{eq:RF} in the form:
\begin{equation*}
E_*(Y)=\HH[k]_*
\end{equation*}
(see \S\ref{subsec:pfRF2}).
The vanishing claim from Theorem \ref{thm:HLCS} follows at once.
Furthermore, the above Rescaling Formula also implies that
\begin{equation*}
\Phi_{2kr}=\dim_{\Q}E_{2kr}(Y)=
\dim_{\Q}\HH_r=\phi_r \, ,
\end{equation*}
for all $r\ge 1$. Replacing $t$ by $t^{2k+1}$ in
the algebraic LCS formula \eqref{eq:klcs} gives
the desired homotopy LCS formula \eqref{eq:hlcs}. 

\subsection{Rationalized loop spaces}
\label{subsec:ratloop}

As is well-known,
the loop space of a simply-connected, finite-type CW--space
has the rational homotopy type of a weak product
of Eilenberg--MacLane spaces of type $K(\Q,n)$.
Thus, Theorem~\ref{thm:HLCS} yields the following.

\begin{cor}
\label{cor:kqn}
Let $Y$ be a finite-type, simply-connected CW--space.
Suppose $Y$ is a homological $k$--rescaling of a finite-type,
connected CW--space $X$, with $H^*(X,\Q)$ Koszul. 
Then:
\begin{equation}
\label{eq:rationalloop}
\Omega Y\simeq_{\Q} \sideset{}{^w}\prod_{r=1}^{\infty} 
K(\Q,2kr)^{\Phi_{2kr}}\, ,
\end{equation}
where $\Phi_{2kr}$ is given by the homotopy LCS formula \eqref{eq:hlcs}.
\end{cor}

Consequently, the rational homotopy type of the loop space of $Y$
is determined by the Poincar\'{e} polynomial of $X$. In particular,
the Poincar\'{e} series of $\Omega Y$ is given by
\begin{equation}
\label{eq:poinseries}
P_{\Omega Y}(t)=P_X(-t^{2k})^{-1}.
\end{equation}
In fact, the Milnor--Moore theorem \cite{MM} insures that
$H_*(\Omega Y,\Q)\cong U(L_*(X)[k])$, as Hopf algebras.

\section{Homological rescalings: existence, uniqueness and simple examples}
\label{sec:surfaces}

In this section, we discuss the existence and uniqueness of
homological rescalings of spaces.  As we shall see, rational homotopy
theory guarantees the existence of such rescalings, but not
their uniqueness, except when certain homological criteria
are satisfied.  We conclude with some simple
cases where the (unique) rescaling can be constructed directly,
via geometric methods. More examples will be given
in \S\S\ref{sec:arrangements} and \ref{sec:links}.

\subsection{Existence and non-uniqueness}
\label{subsec:existence}

Given a connected, finite-type CW--space $X$, and an integer
$k\ge 1$, there exists a simply-connected, finite-type CW--space $Y$
such that $H^*(Y,\Q)$ is the $k$--rescaling of $H^*(X,\Q)$.
Indeed, $(H^*(X,\Q)[k], d=0)$ is a $1$--connected \dga,
with Sullivan minimal model $\mathcal{M}$.  Hence,
there exists a finite-type, simply-connected
CW--space $Y$  such that $\mathcal{M}(Y)=\mathcal{M}$.
In particular, $H^*(Y,\Q)=H^*(X,\Q)[k]$.

By construction, $Y$ is formal.  Hence, $Y$ is uniquely
determined (up to rational homotopy equivalence) among
spaces with the same cohomology ring.  But there may be other,
simply-connected, non-formal spaces, $Z$, with $H^*(Z,\Q)=H^*(Y,\Q)$.
In other words, the homological $k$--rescaling property \eqref{eq:hxy}
alone does not determine the rational homotopy type of $Y$, as the
next example shows.

\begin{example}
\label{ex:nonunique}
Fix an integer $k\ge 1$, and consider the space
$X=S^1\vee S^1 \vee S^{2k+2}$.  Plainly, the formal
$k$--rescaling of $X$ is
$Y=S^{2k+1}\vee S^{2k+1}\vee S^{(2k+1)(2k+2)}$.

Now let
$Z=(S^{2k+1}\vee S^{2k+1})\cup _{\alpha} e^{(2k+1)(2k+2)}$,
where the attaching map of the top cell is the iterated
Whitehead product $\alpha=\ad_x^{2k+2}(y)=[x,[\cdots[x,y]]]$,
and $x,y\in\pi_{2k+1}(S^{2k+1}\vee S^{2k+1})$ are the homotopy
classes of the factors of the wedge.  It is readily seen that
$\pi_1(Z)=0$ and $H^*(Z,\Q)=H^*(X,\Q)[k]$.  Thus, $Z$ is also
a $k$--rescaling of $X$.

The Quillen minimal model of $Z$ is
$(\BL(x,y,z) , \partial)$, with $\deg x=\deg y=2k$, $\deg z=(2k+1)(2k+2)-1$,
and $\partial x=\partial y=0$, $\partial z=\ad _x^{2k+2}(y)$.
Clearly, the differential is not quadratic, and therefore
$Z$ is not formal.
\end{example}

\subsection{Conditions for uniqueness}
\label{subsec:unique}

The above example notwithstanding, there are several commonly
occurring situations where homological rescalings are  unique
(up to rational homotopy equivalence).  We list two such
situations.

\begin{prop}[Shiga-Yagita \cite{SYa}]
\label{sss1}
If $ H^{>d}(X,\Q)=0$, then $X$ has a unique $k$--rescaling,
for all $k>(d-1)/2$.  
\end{prop}

\begin{proof}
The result follows from \cite[Theorem 5.4]{SYa}. (N.B. There is no need to
assume that $2k+1$ is a prime number, as the authors of \cite{SYa} do.)
\end{proof}

\begin{prop}
\label{sss3}
If $H^*(X,\Q)=H^*\big(\prod^{w}_i K(\Q,n_i)\big)$,
then $X$ has a unique $k$--rescaling,
for all $k\ge 1$.
\end{prop}

\begin{proof}
Suppose $Y$ is a simply-connected space which has
the rational cohomology of a (weak) product of Eilenberg-MacLane
spaces, $K=\prod^{w}_i K(\Q,(2k+1) n_i)$.  Then $Y$ admits a classifying
map to $K$.  By the Whitehead-Serre theorem, the map
$Y\to K$ is a rational homotopy equivalence.
\end{proof}

By the discussion in \S\ref{subsec:existence}, once we know $X$
has a unique $k$--rescaling, then such a rescaling, $Y$, must be
a formal space (even if $X$ itself is not formal).

\subsection{Wedges of circles}
\label{subsec:wedges}
Start with $X=S^1$.  This is a formal space, with rationalization
$X_{\Q}=K(\Q,1)$. Hence, $H^*(X,\Q)$ is a Koszul algebra.
For each $k\ge 1$, the (unique up to $\Q$--equivalence)
homological $k$--rescaling of $X$ is $Y=S^{2k+1}$.
Clearly  $E_*(Y)=\BL(x)$, with $\deg x=2k$, and
$\Omega S^{2k+1} \simeq_{\Q} K(\Q,2k)$,
a result that goes back to Serre's thesis.

Next, consider $X=\bigvee^n S^1$.  Note that $\pi_1(X)=\F_n$,
the free group on $n$ generators.  Clearly, $X$ is formal.
Moreover, $X_{\Q}=K(\F_n\otimes \Q,1)$, and thus $H^*(X,\Q)$
is a Koszul algebra.  For each $k\ge 1$,
the (unique) $k$--rescaling of $X$ is $Y=\bigvee^n S^{2k+1}$.
By a well-known result of Magnus, $L_*(X)=\gr^*(\F_n)\otimes \Q$
is the free Lie algebra on $n$ generators. By Theorem~\ref{thm:RF1}:
\begin{equation*}
\label{eq:hliecirc}
E_*(Y)=\BL(x_1,\dots,x_{n}),
\end{equation*}
with $\deg x_i=2k$, which is a particular case of
the classical Hilton--Milnor theorem. The ranks 
of the homotopy Lie algebra of $Y$ may be computed
by means of Theorem~\ref{thm:HLCS}. We find:
$\Phi_q=0$ if $q$ is not divisible by $2k$, and
\begin{equation*}
\label{eq:hlcscirc}
\prod_{r\ge 1} (1-t^{(2k+1)r})^{\Phi_{2kr}}=1-nt^{2k+1}\, .
\end{equation*}
By the Milnor--Moore theorem, 
$H_*(\Omega Y,\Q)=UE_*(Y)=T(x_1,\dots,x_n)$, the
tensor algebra on $x_i$, and
$P_{\Omega Y}(t)=\frac{1}{1-nt^{2k}}$,
a result that goes back to Bott and Samelson \cite{BS}.

\subsection{Products of circles}
\label{subsec:tori}
Another simple situation is that of tori.
Let $X= (S^1)^{\times n}$ be the $n$--torus.
Clearly, $X$ is formal,
and $X_{\Q}=K(\Q^n,1)$; hence $H^*(X,\Q)$ is Koszul.
For each $k\ge 1$, the (unique) $k$--rescaling of $X$ is
 $Y=(S^{2k+1})^{\times n}$.
By Theorem~\ref{thm:RF1}:
\begin{equation*}
\label{eq:hlietorus}
E_*(Y)=\BL^{\ab}( x_1,\dots,x_n):=\BL( x_1,\dots,x_n)/ ([x_i,x_j]=0),
\end{equation*}
the free abelian Lie algebra on generators in degree $2k$
(this also follows
from \cite{CCX}, by noting that $X\simeq M(\mathcal{D}_n)$
and $Y\simeq M(\mathcal{D}_n^{k+1})$,
where $\mathcal{D}_n$ is the Boolean arrangement in $\C^n$).
As above, we compute:
$\Phi_q=0$ if $2k\nmid q$, and
\begin{equation*}
\label{eq:hlcstori}
\prod_{r\ge 1} (1-t^{(2k+1)r})^{\Phi_{2kr}}=(1-t^{2k+1})^n\, .
\end{equation*}
Furthermore,
$\Omega Y \simeq_{\Q} K(\Q,2k)^{\times n}$.
Hence: $H_*(\Omega Y,\Q)=\Q[x_1,\dots ,x_n]$, 
the polynomial algebra on 
$x_i$, and $P_{\Omega Y}(t)=\frac{1}{(1-t^{2k})^n}$.

\subsection{Surfaces}
\label{subsec:surfaces}

Let $X=\#^g S^1\times S^1$ be a compact, orientable
surface of genus $g$.  Being a K\"{a}hler manifold, $X$
is formal, see \cite{DGMS}.
Moreover, $A=H^*(X,\Q)$ is a Koszul algebra
(see for instance \cite[Example~2.2]{BP},
or simply note that $A$ has a quadratic
Gr\"obner basis, which is a well-known sufficient condition
for Koszulness).
For each $k\ge 1$, the (unique) $k$--rescaling
of $X$ is $Y=\#^g S^{2k+1}\times S^{2k+1}$.
Thus, we may apply Theorem~\ref{thm:RF1}
to determine the homotopy Lie algebra of $Y$ from the associated graded
Lie algebra of $\pi_1 X$.
Using a result of Labute \cite{Lb}, we obtain: 
\begin{equation*}
\label{eq:hliesurf}
E_*(Y)=\BL(x_1,\dots,x_{2g})/([x_1,x_2]+\cdots +[x_{2g-1},x_{2g}]=0),
\end{equation*}
with $\deg x_i=2k$.
As above, we compute the ranks of $\pi_q(\Omega Y)$ to be:
$\Phi_q=0$ if $2k \nmid q$, and%
\begin{equation*}
\label{eq:hlcsrs}
\prod_{r\ge 1} (1-t^{(2k+1)r})^{\Phi_{2kr}}=1-2gt^{2k+1}+t^{4k+2}\, .
\end{equation*}
It follows that $P_{\Omega Y}(t)=(1-2gt^{2k}+t^{4k})^{-1}$.

\subsection{Rescaling non-Koszul algebras}
\label{subsec:nonkoszul}
We conclude with a simple example showing that the converse to
Proposition~\ref{prop:coformal} does not hold.
This example also shows that the Koszul condition
is indeed necessary for Theorems \ref{thm:RF2}
and \ref{thm:HLCS} to hold.

\begin{example}
\label{ex:nonkoszul}
Let $X=S^1 \vee S^{2}$, with $k$--rescaling $Y=S^{2k+1}\vee S^{2(2k+1)}$.
Of course, both $X$ and $Y$ are formal.  Moreover, $Y$ is coformal.
Nevertheless, the algebra $A=H^*(X,\Q)$ is not generated
in degree~$1$, and thus, in particular, $A$ is not Koszul.
Finally, note that $\Phi_{4k+1}(Y)=1$, even though $2k\nmid 4k+1$;
in particular, $E_*(Y)\ne L_*(X)[k]$.
\end{example}

\part{\Large\bf The Malcev formula}
\label{part:MF}

In this second part, we study  groups of based
homotopy classes of maps between loop spaces, and prove
Theorem~\ref{thm:MF} from the Introduction, using techniques
of Lazard, Quillen, and Baues.

The goal is to establish filtered group isomorphisms between
 the Malcev completion of $\pi_1(X)$,
the completion of $[\Omega S^{2k+1},\Omega Y]$,
and the Milnor--Moore group of coalgebra maps from
$H_*(\Omega S^{2k+1},\Q)$ to $H_*(\Omega Y,\Q)$,
all under the assumption that $Y$ is a homological
$k$--rescaling of a {\em formal}\/ space $X$ with
Koszul cohomology algebra.

The Malcev Formula,
$\pi_1(X)\otimes \Q\cong \comp{[\Omega S^{2k+1},\Omega Y]}$,
upgrades the Rescaling Formula to the level of filtered groups.
Remarkably, this formula is strong enough to insure the
formality of $X$.  A parallel analysis gives information 
on the ``group of homotopy groups",
$[\Omega S^{2},\Omega Y]$, of Cohen and Gitler.

\section{Malcev completions and Baues formula}
\label{sec:malcev}

This section contains our main technical tool for the proof
of Theorem~\ref{thm:MF}:
a formula (due to Baues \cite{Ba}), identifying, under certain
conditions, the Malcev completion of the group of
based homotopy classes $[K,\Omega Y]$ with an
algebraically defined exponential group.

\subsection{Complete Lie algebras and Campbell--Hausdorff groups}
\label{subsec:ch}

Our results involve {\em Malcev completions}
of groups.  We will need the Lie algebra form of Malcev completion,
so we start by reviewing the relevant material from
\cite[Appendix~A]{Q}  (see also \cite{La, Se}).

By a {\em complete Lie algebra} (\clie) we mean an (ungraded!)
Lie algebra $L$, together with a complete, descending $\Q$--vector
space filtration, $\{F_rL\}_{r \ge 1}$, satisfying $F_1L=L$ and
\begin{equation}
\label{eq:clie}
[F_1L,F_rL]\subset F_{r+1}L, \quad \text{for all $r$}.
\end{equation}
Completeness of the filtration $\{F_rL\}_{r \ge 1}$  means that the induced 
topology on $L$ is Hausdorff, and that every Cauchy sequence converges.
In other words, the canonical map to the inverse limit,
$\pi\co L\to \varprojlim_r L/F_r L$, is a vector space isomorphism.

Let $L$ and $L'$ be two \clie's.
We say that a Lie algebra map, $f\co L\to L'$, is a weak \clie~morphism
if $f(F_r L)\subset F_rL'$, for $r$ sufficiently large.
If $f(F_r L)\subset F_rL'$ for all $r$, we simply say that $f$ is a
\clie~morphism.

To a complete Lie algebra $L$, one associates, in a functorial way,
a filtered group, called the {\em exponential group of $L$}, and
denoted by $\exp(L)$. The underlying set of $\exp(L)$ is just $L$,
while the group law is given by the classical
Campbell--Hausdorff multiplication from local Lie theory:
\begin{equation}
\label{eq:chmult}
x \cdot y=x+y+\tfrac{1}{2}[x,y]+
\tfrac{1}{12}[x,[x,y]]+\tfrac{1}{12}[y,[y,x]]+\cdots ,
\quad \text{for $x,y\in L$}.
\end{equation}
(The convergence of the series follows from condition \eqref{eq:clie},
together with the completeness of the filtration topology.)
The Lie filtration, $\{F_rL\}_{r \ge 1}$, passes to a filtration by
normal subgroups, $\{\exp(F_rL)\}_{r \ge 1}$, of $\exp (L)$.

\begin{remark}
\label{rem:cvsm}
In \cite{Q}, Quillen uses a special kind of \clie's, the so-called
{\em Malcev Lie algebras} (\mlie's).  In his definition,
condition \eqref{eq:clie} above is replaced by the stronger condition
\begin{equation*}
\label{eq:mlie1}
[F_rL,F_sL]\subset F_{r+s}L, \quad \text{for all $r$ and $s$},
\end{equation*}
which implies that $\gr^*(L):=\bigoplus _{r\ge 1} F_rL/F_{r+1}L $
has a natural \grlie~structure.  Additionally, Quillen assumes that
the Lie algebra $\gr^*(L)$ is generated by $\gr^1(L)$.  For example,
any nilpotent Lie algebra $L$, with lower central series filtration
$\{\Gamma_r L\}_{r \ge 1}$, is a Malcev Lie algebra.
We have chosen to isolate, in our definition, the minimal conditions
needed for the construction of the exponential group.
\end{remark}

\subsection{Malcev completion of groups}
\label{subsec:malcev}
The Malcev completion functor of Quillen \cite{Q} sends
a group $G$ to the filtered group $G\otimes \Q:=\exp (M_G)$,
where $M_G$ is the {\em Malcev Lie algebra of $G$}, constructed
as follows.

Let $\Q G$ be the rational group algebra of $G$, endowed
with the $I$--adic filtration $\{I^r \Q G\}_{r\ge 0}$, where
$I$ is the augmentation ideal.  The inverse limit
$\widehat{\Q G}=\varprojlim_r \Q G/I^r \Q G$ has a natural complete
Hopf algebra structure.  By definition, $M _G$ is the Lie algebra of
primitive elements in $\widehat{\Q G}$.
If $G$ is a nilpotent group, then $M_G$ is a nilpotent Lie algebra,
and $G\otimes \Q$ coincides with the usual Malcev completion of $G$.
In general, $G\otimes \Q=\varprojlim_n ( (G/\Gamma_n G) \otimes \Q )$.

\begin{example}
\label{ex:mformal}
Let $\HH_*$ be a \grlie~, generated by $\HH_1$.
Denote by $\compp{\HH}$ the completion of $\HH_*$ with respect to
the degree filtration.  The elements of $\compp{\HH}$ are formal
series of the form
\begin{equation}
\label{eq:f1}
s=\sum_{i\ge 1} s_i,\quad \text{with $s_i\in \HH_i$}.
\end{equation}
The Lie bracket on $\HH_*$ extends bilinearly to the completion,
making $\compp{\HH}$ into a Malcev Lie algebra, with filtration
given by formal series order:
\begin{equation}
\label{eq:f2}
F_r\compp{\HH}:=\{s\mid s_i=0\ \text{for $i<r$}\}.
\end{equation}

Let $X$ be a connected CW--space of finite type, with cohomology
algebra $A^*=H^*(X,\Q)$, and holonomy Lie algebra $\HH_*=\HH_*(A)$.
If $X$ is a formal space, then 
\begin{equation}
\label{eq:f3}
\pi_1(X)\otimes \Q=\exp(\compp{\HH}_*)=\exp (\comp{L_*(X)}) ,
\end{equation}
see for instance \cite[Lemma~1.8]{MP}. 
\end{example}

\subsection{The exponential formula of Baues}
\label{subsec:baues}

To state Baues' formula,  we start
by recalling the setup from \cite[\S{VI.1}]{Ba}.

Let $C_*$ be a connected, cocommutative, graded $\Q$--coalgebra
(\cgc) of finite type. That is, $C_*$ is the dual of a finite-type \cga,
with diagonal $\Delta\co C_*\to C_*\otimes C_*$ dual to the
multiplication. Set $C_+=\bigoplus _{k>0} C_k$, and let
$\overline{\Delta}\co C_+\to C_+\otimes C_+$ be the reduced
diagonal.

Let $E_*$ be a \glie~of finite type, with Lie bracket denoted by
$b\co E_*\otimes E_*\to E_*$. Denote by $\Hom(C_+,E_*)$ the
$\Q$--vector space of degree~$0$ linear maps
from $C_+$ to $E_*$.  For $f, g\in \Hom(C_+,E_*)$, define
$[f,g]\in \Hom(C_+,E_*)$ to be the composite
\begin{equation}
\label{eq:hombr}
 C_+\xrightarrow{\:\overline{\Delta}\:} C_+\otimes C_+
\xrightarrow{f\otimes g} E_*\otimes E_*
\xrightarrow{\: b\: } E_*.
\end{equation}
Endowed with this bracket, $\Hom(C_+,E_*)$ becomes
an (ungraded) Lie algebra. Moreover, if $\dim_{\Q} C_*<\infty$, then
$\Hom(C_+,E_*)$ is a nilpotent Lie algebra, so we may speak of the
exponential group $\exp(\Hom(C_+,E_*))$.
\begin{theorem}[Baues~\cite{Ba}, Theorem~VI.1.3]
\label{thm:baues}
Let $K$ be a connected, finite complex, with homology coalgebra $H_*(K,\Q)$,
dual to the algebra $H^*(K,\Q)$.
Let $Y$ be a simply-connected CW--space
of finite type, with rational homotopy Lie algebra $E_*(Y)=\pi_*(\Omega Y)\otimes \Q$.
There is then a group isomorphism
\[
[K,\Omega Y] \otimes \Q \cong \exp (\Hom (H_+(K,\Q),E_*(Y))),
\]
natural in both $K$ and $Y$.
\end{theorem}

\section{Groups of homotopy classes and the Milnor-Moore groups}
\label{sec:thmD}

In this section, we prove Theorem~\ref{thm:MF} from the Introduction.
Let $Y$ be a based, simply-connected CW--space of finite type.
Our first aim is to establish
a precise relationship between the group
$\comp{[\Omega S^{m}, \Omega Y]}$ and the exponential
group of $\comp{E_*(Y)\{m\}}$ (the completion of the
``rebracketed" homotopy Lie algebra of $Y$).  From this,
half of Theorem~\ref{thm:MF} will follow.

Our second aim is to establish  an isomorphism between filtered groups
of the form $\comp{[K,\Omega Y]}$, and the Milnor--Moore groups
$\Hom^{\coalg} (H_*(K,\Q),H_*(\Omega Y,\Q))$. This will finish
the proof of Theorem~\ref{thm:MF}.

\subsection{Rebracketing the homotopy Lie algebra}
\label{subsec:mres}

Let $E_*$ be a graded Lie algebra.  For a fixed integer $m\ge 2$,
consider the Lie subalgebra $E_*\{m\}=\bigoplus _{r\ge 1} E_{r(m-1)}$.
Now modify the Lie bracket on $E_*\{m\}$ by making all the elements
of odd degree commute, while leaving the other brackets unchanged.
This modified bracket turns $E_*\{m\}$ into a \grlie~(with the usual
signs in the Lie identities). The completion, $\comp{E_*\{m\}}$,
consists of elements of the form
\begin{equation}
\label{eq:e1}
s=\sum_{i\ge 1} s_i,\quad \text{with $s_i\in E_{i(m-1)}$},
\end{equation}
with Lie bracket induced from the \grlie~bracket on $E_*\{m\}$.
Endowed with the filtration
\begin{equation}
\label{eq:e2}
F_r\comp{E_*\{m\}}:=\{s\mid s_i=0\ \text{for $i<r$}\},
\end{equation}
$\comp{E_*\{m\}}$ becomes a \clie.

\begin{theorem}
\label{thm:exp}
Let $Y$ be a based, simply-connected CW--space of finite type,
with homotopy Lie algebra $E_*=E_*(Y)$. Then, for each $m\ge 2$,
there is a filtered group isomorphism,
\[
\comp{[\Omega S^m, \Omega Y]} \cong \exp (\comp{E_*\{m\}}\,).
\]
\end{theorem}

\subsection{Proof of Theorem~\ref{thm:exp}}
\label{subsec:thmexp}

By a well-known result of Milnor--Moore~\cite{MM}, the coalgebra
$C_*=H_*(\Omega S^m,\Q)$ is, in fact, a Hopf algebra, freely generated
(as an algebra) by a primitive element, $v\in C_{m-1}$.
Consequently, we have the following simple description of the
coalgebra structure on $C_*$.  If $m$ is odd, then:
\begin{equation}
\label{eq:odd}
\Delta(v^k)=\sum_{r=0}^{k} \binom{k}{r} v^r \otimes v^{k-r},
\quad \text{for all $k$}.
\end{equation}
If $m$ is even, then:
\begin{equation}
\label{eq:even}
\begin{cases}
\Delta(v^{2k})&=\sum_{r=0}^{k} \tbinom{k}{r} v^{2r} \otimes v^{2(k-r)},
\\[2pt]
\Delta(v^{2k+1})&=\sum_{r=0}^{k} \tbinom{k}{r}
(v^{2r+1} \otimes v^{2(k-r)} + v^{2r} \otimes v^{2(k-r)+1}).
\end{cases}
\end{equation}

Now filter the complex $K=\Omega S^m$,
by choosing $K_r$ to be the $r(m-1)$-th skeleton of $K$. 
It is readily seen that
\begin{equation}
\label{eq:sfilt}
H_*(K_r,\Q)=\Q\text{-span} \{ v^i \mid i\le r \},
\quad \text{for all $r\ge 0$}.
\end{equation}

The inclusion $K_{r-1}\to K$ induces a homomorphism 
$\gamma_r\co \Hom(H_+(K),E_*)\to \Hom(H_+(K_{r-1}),E_*)$,
which obviously  preserves Lie brackets\footnote{Recall
we are using $\Q$--coefficients, unless otherwise specified.}.
Since the filtration of $K$ is exhaustive, we have
$H_*(K)=\varinjlim_r H_*(K_r)$. Hence, the
maps $\{\gamma_r\}_{r\ge 1}$ define an isomorphism of Lie algebras,
\begin{equation}
\label{eq:lieiso}
\gamma\co \Hom (H_+(K),E_*(Y))\longrightarrow
\varprojlim_r \Hom (H_+(K_{r-1}),E_*(Y)).
\end{equation}

The next three Lemmas finish the proof of Theorem~\ref{thm:exp}.

\begin{lemma}
\label{lem:exp1}
Endowed with the inverse limit filtration $\{F_r\}_{r\ge 1}$
coming from the isomorphism
\eqref{eq:lieiso},
$\Hom(C_+,E_*)$ becomes a complete Lie algebra.
\end{lemma}

\begin{proof}
From the way it was defined, $\{F_r\}_{r\ge 1}$ is a complete
filtration. Since $K_0=\text{point}$, $F_1=\Hom(C_+,E_*)$.
It remains to check condition \eqref{eq:clie}, using 
definition \eqref{eq:hombr} of the Lie bracket.

Equations \eqref{eq:odd}--\eqref{eq:sfilt} together imply that
\begin{equation}
\label{eq:Qufilt}
\overline{\Delta} (H_+(K_r)) \subset H_+(K_{r-1})\otimes H_+(K_{r-1}),
\quad \text{for all $r\ge 1$}.
\end{equation}
Condition \eqref{eq:clie} now follows from \eqref{eq:Qufilt}.
\end{proof}

\begin{lemma}
\label{lem:exp2}
There is an isomorphism of filtered groups,
\[\comp{[\Omega S^m,\Omega Y]} \cong \exp (\Hom(C_+,E_*)).\]
\end{lemma}

\begin{proof}
Recall that $K=\Omega S^m$,  $C_*=H_*(K)$, and  $E_*=E_*(Y)$.
Baues' formula (Theorem \ref{thm:baues}) implies that
the filtered group $\comp{[K,\Omega Y]}$ is isomorphic to
$\varprojlim_r \exp (\Hom(H_+(K_{r-1}),E_*))$,
both groups being endowed with the canonical
inverse limit filtration.

By Lemma~\ref{lem:exp1}, the exponential group
$\exp(\Hom(C_+,E_*))$ is defined.  From the definition
of the filtration $\{F_r\}_{r\ge 1}$ on
$\Hom(C_+,E_*)$, the map $\gamma_r$ sends $F_r$ to $0$.
Hence, $\gamma_r$ eventually preserves filtrations,
and thus is a weak \clie~morphism. Consequently, $\gamma_r$ 
induces a map of exponential groups,
$\exp(\gamma_r)\co \exp(\Hom(C_+,$ $E_*))\to 
\exp (\Hom(H_+(K_{r-1}),E_*))$.
Passing to the inverse limit, we obtain a group isomorphism
\[
\exp(\gamma)\co \exp(\Hom(C_+,E_*))\xrightarrow{\:\simeq\:}
\varprojlim_r \exp (\Hom(H_+(K_{r-1}),E_*)).
\]
By definition
of the respective filtrations, the map $\exp(\gamma)$ is 
filtration-preser\-ving.
\end{proof}

\begin{lemma}
\label{lem:exp3}
The filtered Lie algebra $\Hom(C_+,E_*)$ is isomorphic to $\comp{E_*\{m\}}$.
\end{lemma}

\begin{proof}
Let $\{c_k\}_{k\ge 1}$ be an arbitrary sequence of non-zero rational numbers.
The map sending $f\in \Hom(C_+,E_*)$ to the formal series
$s=\sum_{k\ge 1} c_kf(v^k)\in \comp{E_*\{m\}}$ is an isomorphism
between the underlying vector spaces.
Using \eqref{eq:sfilt}, we see that
the respective filtrations are preserved under this map.
Finally, inspection of equations \eqref{eq:odd}--\eqref{eq:even}
shows that the Lie brackets are preserved as well, provided we take
$c_k=\frac{1}{k!}$ (when $m$ is odd), or
$c_{2k}=c_{2k+1}=\frac{1}{k!}$ (when $m$ is even).
\end{proof}

\subsection{Proof of Theorem~\ref{thm:MF}}
\label{subsec:pfMF}

We are now in position to prove Theorem~\ref{thm:MF} from
the Introduction.

\begin{prop}
\label{prop:D1}
Let $X$ and $Y$ be finite-type CW--spaces ($X$ connected and
$Y$ simply-connected), such that $H^*(Y,\Q)=H^*(X,\Q)[k]$,
for some $k\ge 1$. Assume that $H^*(X,\Q)$ is a Koszul algebra.
Then there is an isomorphism of filtered groups,
\[
\comp{[\Omega S^{2k+1}, \Omega Y]}\cong \exp (\comp{L_*(X)}).
\]
\end{prop}

\begin{proof}
By Theorem~\ref{thm:exp},
it is enough to verify that $\comp{E_*\{2k+1\}}\cong \comp{L_*(X)}$, 
as filtered Lie algebras. 
This follows directly from Theorem~\ref{thm:RF2}.
The main point is that the Rescaling Formula  
$E_*(Y)=L_*(X)[k]$ holds, thus
implying that the constructions \eqref{eq:e1}--\eqref{eq:e2} and
\eqref{eq:f1}--\eqref{eq:f2} define isomorphic filtered Lie algebras.
\end{proof}

\begin{cor}
\label{cor:Mform}
Let $X$ and $Y$ satisfy the assumptions 
from Proposition~\ref{prop:D1}. Then
there is an isomorphism of filtered groups,
\begin{equation}
\label{eq:MalF}
\pi_1(X)\otimes \Q \cong \comp{[\Omega S^{2k+1}, \Omega Y]} \, ,
\end{equation}
if and only if $X$ is formal.
\end{cor}

\begin{proof}
In one direction, the formality of $X$ gives a filtered group isomorphism,
$\pi_1(X) \otimes \Q \cong \exp (\comp{L_*(X)})$; see \eqref{eq:f3}. We may
conclude by applying the above Proposition.

Conversely, assume that the Malcev Formula \eqref{eq:MalF} holds.
By Proposition \ref{prop:D1}, this implies that
$\pi_1(X) \otimes \Q \cong \exp (\comp{L_*(X)})$.
Using Lemma~\ref{lem:r1}, we infer that
$\pi_1(X) \otimes \Q \cong \exp (\compp{\HH}_*)$, where $\HH_*$ is the 
holonomy Lie algebra of $A^* = H^*(X, \Q)$. In terms of
Sullivan models, this isomorphism translates to the fact that $X$ has 
the same $1$--minimal model as the \dga~ $(A^*,0)$. On the other hand, 
one knows from \cite[Proposition~5.2]{PY} that both the minimal model of $X$ 
and that of $(A^*,0)$ are generated in degree~$1$. The formality of 
$X$ follows at once.
\end{proof}

The Corollary proves the implication \eqref{Da} $\Longrightarrow$ 
\eqref{Db} from Theorem~\ref{thm:MF}, and half of the other 
implication.  The next Proposition ends the proof
of Theorem~\ref{thm:MF}, once we apply it to
the loop space $K=\Omega S^{2k+1}$, endowed with the
Morse--James filtration.

\begin{prop}
\label{prop:homco}
Let $K$ be a based CW--complex, endowed with an exhaustive filtration
by connected, finite subcomplexes, $\{K_r\}_{r\ge 0}$,
with $K_0=*$, and let $Y$ be a based, simply-connected
CW--space of finite type.  Then there is an isomorphism of
filtered groups,
\[
\comp{[K,\Omega Y]}\cong \Hom^{\coalg} (H_*(K,\Q),H_*(\Omega Y,\Q)).
\]
\end{prop}

\begin{proof}
From the definition of the filtration~\eqref{eq:filt1}, we see that
the Milnor--Moore group of coalgebra maps from 
$H_*(K,\Q)$ to $H_*(\Omega Y,\Q)$ is isomorphic to the group 
$\varprojlim_r \Hom^{\coalg}  (H_*(K_{r-1},\Q),H_*(\Omega Y,\Q))$, 
filtered as an inverse limit. 

Now recall that
$\comp{[K,\Omega Y]}\cong\varprojlim_r ([K_{r-1},\Omega Y]\otimes \Q)$.
Thus, it is enough to show that the ``Hurewicz" homomorphism,
$h\co [K_{r-1},\Omega Y]\to   
\Hom^{\coalg}(H_*(K_{r-1}),$ $H_*(\Omega Y))$,
defined by $h(f)=f_*$, gives rise to a natural isomorphism
\[
h_{\Q}\co [K_{r-1},\Omega Y]\otimes \Q \longrightarrow  
  \Hom^{\coalg}  (H_*(K_{r-1},\Q),H_*(\Omega Y,\Q)).
\]

Since $K_{r-1}$ is a finite complex, Theorem~II.3.11 in 
Hilton--Mislin--Roitberg \cite{HMR} gives a group isomorphism
\begin{equation*}
\label{hmr}
[K_{r-1},\Omega Y]\otimes \Q \cong [K_{r-1}, \Omega Y_{\Q}].
\end{equation*}
Thus, we may use Proposition~1 in Scheerer \cite{Sc} to conclude that
the Hurewicz map $h_{\Q}$ is a group isomorphism.
\end{proof}

\section{Groups of homotopy groups}  
\label{sec:glorified}

We now apply the results from \S\ref{sec:thmD} to the
study of the ``group of homotopy groups" of a loop space,
$[\Omega S^2,\Omega Y]$, and its close approximation,
the Milnor--Moore group of coalgebra maps,
$\Hom^{\coalg} (H_*(\Omega S^2,\Q),H_*(\Omega Y,\Q))$.

\subsection{Loop space of $S^2$}
\label{subsec:loopsphere}

We start with a situation considered
by Cohen and Gitler~\cite{CG2},
in the context of configurations spaces.

\begin{prop}
\label{prop:answer} 
Let $X$ and $Y$ be finite-type CW--spaces ($X$ connected and
$Y$ simply-connected), such that $H^*(Y,\Q)=H^*(X,\Q)[k]$,
for some $k\ge 1$. Assume that $H^*(X,\Q)$ is a Koszul algebra.
If $X$ is formal, then there is a group isomorphism
\begin{equation}
\label{gpiso}
\Hom^{\coalg} (H_*(\Omega S^2,\Q),H_*(\Omega Y,\Q))\cong \pi_1(X)\otimes \Q.
\end{equation}
\end{prop}

\begin{proof}
By Proposition~\ref{prop:homco} and Theorem~\ref{thm:exp},
the group on the left side is isomorphic to
$\comp{[\Omega S^2,\Omega Y]}=\exp (\comp{E_*\{2\}})$, where $E_*=E_*(Y)$.
By Theorem~\ref{thm:RF2} and formula \eqref{eq:second}, there
is an isomorphism of Lie algebras, $\comp{E_*\{2\}}\cong \compp{\HH}_*$,
where $\HH_*$ is the holonomy Lie algebra of $X$.
Even though this isomorphism does not preserve the
filtrations on $\comp{E_*\{2\}}$ and $\compp{\HH}_*$,
it does preserve the corresponding topologies.
Thus, it induces an isomorphism between $\exp(\comp{E_*\{2\}})$
and $\exp(\compp{\HH}_*)=\pi_1(X)\otimes \Q$.
\end{proof}

\begin{remark}
\label{rem:CGS1}
This Proposition answers a question posed by Cohen and Gitler 
in \cite[\S6]{CG2}, in the particular case when $X=M(\B_{\ell})$ 
is the complement of the braid arrangement in $\C^{\ell}$, and 
$Y=M(B_{\ell}^{k+1})$.  As noted in the proof, the isomorphism 
\eqref{gpiso} is not filtration-preserving.  On the other
hand, if one replaces $S^2$ by $S^{2k+1}$ this drawback 
disappears, see Theorem~\ref{thm:MF}.
\end{remark}

\begin{example}
\label{ex:malsimp}
When $Y$ is a wedge of $n$ copies of $S^{2k+1}$ as 
in \S\ref{subsec:wedges}, we get:
\[
\Hom^{\coalg}(H_*(\Omega S^2,\Q),H_*(\Omega 
\big(\bigvee\nolimits^n S^{2k+1}\big),\Q))=
\F_n\otimes \Q.
\] 
This recovers a result of Sato~\cite{Sa}  
(see also \cite[Theorem~5.3]{CG2} and \cite[Theorem~2.7]{CS}).

When $Y$ is a product of $n$ copies of $S^{2k+1}$ as in \S\ref{subsec:tori},
we get:
\[
\Hom^{\coalg}(H_*(\Omega S^2,\Q),H_*(\Omega ((S^{2k+1})^{\times n}),\Q))=\Q^n.
\]
\end{example}

\subsection{Associated graded Lie algebras of exponential groups}
\label{subsec:assoc}
Let $L$ be a complete Lie algebra, with filtration $\{F_r\}_{r\ge 1}$
satisfying the stronger condition from Remark~\ref{rem:cvsm},
\begin{equation}
\label{eq:mlie}
[F_rL,F_sL]\subset F_{r+s}L, \quad \text{for all $r$ and $s$}.
\end{equation}
One knows from Lazard \cite{La}
that the exponential  group $\exp(L)$,
filtered by  the normal subgroups $\{\exp(F_r)\}_{r\ge 1}$,
determines the filtered Lie algebra $L$.

Condition \eqref{eq:mlie} implies that the associated graded,
$\gr^*_F(L):=\bigoplus_{r\ge 1} F_r/F_{r+1}$,
is a \grlie, with Lie bracket induced from $L$. 
It also implies that
\begin{equation}
\label{eq:expf}
(\exp(F_r), \exp(F_s))\subset \exp(F_{r+s}), \quad
\text{for all $r$ and $s$}.
\end{equation}
Hence, the associated graded,
$\gr^*_F(\exp(L))=\bigoplus_{r\ge 1} \exp(F_r)/\exp(F_{r+1})$
is a also a \grlie, with Lie bracket induced from the Campbell--Hausdorff
group commutator.  Moreover, the two associated \grlie's
are (functorially) isomorphic:
\begin{equation}
\label{eq:lazard}
\gr^*_F(L)\cong \gr^*_F(\exp(L)).
\end{equation}
See \cite{La} for full details.

\begin{example}
\label{ex:grlie}
Let $E_*$ be a graded Lie algebra, and let $\comp{E_*\{m\}}$
be the \clie~constructed in \S\ref{subsec:mres}.
The filtration \eqref{eq:e2} of $\comp{E_*\{m\}}$ clearly
satisfies condition \eqref{eq:mlie}.
\end{example}

\begin{remark}
\label{rem:upgrade}

Let $X$ and $Y$ be as in Theorem~\ref{thm:MF}, with $X$ formal.  
The filtered group isomorphism
\[
\comp{[\Omega S^{2k+1}, \Omega Y]}\cong
\pi_1(X) \otimes \Q
\]
established in that theorem passes to the associated graded, 
giving the \glie~isomorphism
\begin{equation}
\label{eq:rescagain}
E_*(Y)\{2k+1\}\cong L_*(X)[k].
\end{equation}

Indeed, we know from Theorem~\ref{thm:exp}
that the filtered group $\comp{[\Omega S^{2k+1}, \Omega Y]}$ is
isomorphic to $\exp(\comp{E_*(Y)\{2k+1\}})$. Moreover, by 
Example~\ref{ex:grlie}, the complete Lie algebra $\comp{E_*(Y)\{2k+1\}}$ 
satisfies condition \eqref{eq:mlie}.  Likewise, we know from \eqref{eq:f3} 
that $\pi_1(X)\otimes \Q=\exp (\comp{L_*(X)})$.
Obviously, the \mlie~$\comp{L_*(X)}$ also satisfies \eqref{eq:mlie}.
From the Lazard formula \eqref{eq:lazard}, we infer that the
\grlie's $\gr^*_F(\comp{E_*(Y)\{2k+1\}})=\bigoplus _{r\ge 1} E_{2kr}(Y)$ 
and $\gr^*_F(\comp{L_*(X)})=\bigoplus _{r\ge 1} L_r(X)$  are isomorphic.  
Formula \eqref{eq:rescagain}  follows.

From the vanishing claim in Theorem~\ref{thm:HLCS},
we know that $E_{*}(Y)\{2k+1\}=E_{*}(Y)$. 
Hence, formula \eqref{eq:rescagain} is actually 
equivalent to the Rescaling Formula \eqref{eq:RF}.
\end{remark}

\subsection{The loop space of $S^2$, again}
\label{subsec:cgs}
Let $Y$ be a simply-connected, finite-type CW--space.
From Proposition~\ref{prop:homco} and Theorem~\ref{thm:exp}, 
we know that the group
$\Hom^{\coalg} (H_*(\Omega S^2,\Q),H_*(\Omega Y,\Q))$, endowed with
the Cohen--Gitler filtration \eqref{eq:filt1} is isomorphic to
$\exp (\comp{E_*(Y)\{2\}})$, where $\comp{E_*(Y)\{2\}}$
is filtered as in \eqref{eq:e2}. By  Example~\ref{ex:grlie},
this filtration satisfies condition \eqref{eq:mlie}.

\begin{cor}
\label{grglor}
Let $Y$ be a simply-connected, finite-type CW--space. Then:

\begin{enumerate}
\item \label{1}
There is an isomorphism of \grlie's,
\[
\gr^*_F(\Hom^{\coalg} (H_*(\Omega S^2,\Q),H_*(\Omega Y,\Q)))\cong E_*(Y)\{2\}.
\]
\item \label{2}
If $Y$ is a $k$--rescaling of a finite-type CW--space $X$,
with $H^*(X,\Q)$ Koszul, then:
\[
 E_*(Y)\{2\}=E_*(Y)=L_*(X)[k].
\]
\end{enumerate}
\end{cor}

\begin{proof}
\eqref{1} By \eqref{eq:lazard}, and the discussion 
above, the associated graded Lie algebra 
$\gr^*_F(\Hom^{\coalg} (H_*(\Omega S^2,\Q), H_*(\Omega
Y, \Q)))$  is isomorphic to the graded Lie algebra 
$\gr^*_F (\comp{E_*(Y)\{2\}})=E_*(Y)\{2\}$.

\eqref{2}  Follows from the definition of rebracketing 
and Theorem~\ref{thm:RF2}.
\end{proof}

\begin{remark}
\label{rem:CGS2}
From Part~\eqref{1} of the Corollary, we recover a result of
Cohen and Sato \cite{CS} (see also \cite[Theorem 5.1, Parts (3) 
and (8)]{CG2}), at least over the rationals.

When $X=M(\B_{\ell})$ and $Y=M(\B_{\ell}^{k+1})$, we recover
from Parts~\eqref{1} and \eqref{2} another result from
\cite{CS} (see also \cite[Theorem 6.1)]{CG2}):
\[
\gr^*_F(\Hom^{\coalg} (H_*(\Omega S^2,\Q),
H_*(\Omega M(\B_{\ell}^{k+1}),\Q))) \cong
\gr^*(\pi_1 M(\B_{\ell}))\otimes \Q [k]\, .
\]
\end{remark}

\begin{remark}
\label{rem:filt}
In general, the filtered groups
\[
\comp{[\Omega S^2, \Omega Y]}=\exp(\comp{E_*(Y)\{2\}}\,)
\quad \text{and}\quad
[\Omega S^2, \Omega Y]\otimes \Q=\exp(M_{[\Omega S^2, \Omega Y]})
\]
do not coincide. (Here recall that $M_G$ stands for the 
Malcev Lie algebra of a group $G$.)

For example, take $Y=S^{2k+1}$, $k\ge 1$.  Suppose
$\comp{[\Omega S^2, \Omega Y]}=\exp(M)$, where $M$ is an \mlie.
From \eqref{eq:lazard}, we infer that the \grlie's $E_*(Y)\{2\}$
and $\gr^*(M)$ are isomorphic.  In particular, $E_*(Y)\{2\}$
is generated by $E_1(Y)$, see Remark~\ref{rem:cvsm}. But $E_1(Y)=0$,
since $E_*(Y)$ is a free Lie algebra on a generator of degree $2k$. 
\end{remark}

\part{\Large\bf Arrangements and links}
\label{part:arrlk}

In this last part, we study homological rescalings in two
particularly interesting geometrical settings:
complements of complex hyperplane arrangements
and complements of classical links. Both classes
of spaces admit naturally defined rescalings.

The Rescaling Formula holds for Koszul (in particular, for
supersolvable) arrangements, as well as for links
with connected linking graph.  On the other hand,
for generic arrangements the Rescaling Formula fails
(due to the non-coformality of the rescaled complement,
detected by higher-order Whitehead products),
while for links, the Malcev Formula can fail (due to the non-formality
of the link complement, detected by Campbell--Hausdorff invariants),
even when the Rescaling Formula does hold.

\section{Rescaling hyperplane arrangements}
\label{sec:arrangements}

In this section, we discuss the Rescaling Formula in the case of
complements of complex hyperplane arrangements with Koszul 
cohomology algebras, as well as the reasons why this formula 
fails for non-Koszul arrangements.

Let $\A$ be a hyperplane arrangement in $\C^{\ell}$, with complement
$X=M(\A)$.  As explained in \S\ref{subsec:CCX}, a homological 
$k$--rescaling of $X$  is provided by $Y=M(\A^{k+1})$, where 
$\A^{k+1}$ is the corresponding redundant subspace arrangement 
in $\C^{(k+1)\ell}$. (By Proposition~\ref{sss1},
this rescaling is unique up to $\Q$--equivalence, for $2k+1>\ell$.)
Recall also from \S\ref{subsec:formal} that both $X$ and $Y$ 
are formal spaces.

\subsection{Koszul arrangements}
\label{subsec:koszul}
By Theorem~\ref{thm:RF1}, the Rescaling Formula,
$E_*(Y)=L_*(X)[k]$, applies precisely to the class of arrangements
for which $H^*(X,\Q)$ is a Koszul algebra. In fact,
by Theorem~\ref{thm:MF}, the stronger Malcev Formula,
$\pi_1(X)\otimes\Q =\comp{[\Omega S^{2k+1}, \Omega Y]}$,
also applies in this case.

Presently, the only arrangements which are known to be Koszul
are the fiber-type (or, supersolvable)
arrangements, cf~\cite{SY}. For such arrangements, the
Rescaling Formula was first established in \cite{CCX}.
Koszulness is known to be equivalent to supersolvability,
within the classes of hypersolvable
and graphic arrangements  (see \cite{JP} and \cite{SS}, respectively).
In general, though, it remains an open question whether there exists a
non-supersolvable arrangement with Koszul cohomology algebra,
see \cite[Problem 6.7.1]{Y01}.

By Proposition~\ref{prop:coformal}, if $H^*(X,\Q)$ is Koszul,
then $Y$ is coformal.  We know from Example~\ref{ex:nonkoszul} that
the converse is not true in general:  a space $X$ may
have a coformal rescaling, $Y$, even though the algebra
$H^*(X,\Q)$ is not generated in degree~$1$ (and thus not Koszul).  
But examples of this sort cannot occur among arrangement complements, 
since their cohomology algebras are always generated in degree~$1$.  
This leads us to pose the following.

\begin{question}
\label{quest:cofkos}
Let $X=M(\A)$ and $Y=M(\A^{k+1})$.  If $Y$ is coformal, 
is $H^*(X,\Q)$ Koszul?
\end{question}

Now let $\A$ be a fiber-type arrangement, with exponents
$d_1,\dots ,d_{\ell}$.   As is well-known, the Poincar\'e polynomial
of $X$ factors as $P_X(t)=\prod_{i=1}^{\ell} (1+d_i t)$, see \cite{FR}.
Let $\Phi_n=\rank \pi_n(\Omega Y)$.  From Theorem~\ref{thm:HLCS}, 
we obtain: $\Phi_q=0$ if $2k\nmid q$, and
\[
\prod_{r\ge 1} (1-t^{(2k+1)r})^{\Phi_{2kr}}=
\prod_{i=1}^{\ell} (1-d_i t^{2k+1}).
\]
Therefore, the ranks of $\pi_*(Y)$---and, hence, the
rational  homotopy type of $\Omega Y$---are determined
solely by the exponents of $\A$. Furthermore,
\[
P_{\Omega Y}(t)=\prod_{i=1}^{\ell} (1-d_i t^{2k})^{-1}.
\]

\subsection{Generic arrangements}
\label{subsec:generic}
We now illustrate the failure of the Rescaling Formula 
in the case of arrangements of hyperplanes in general position.
For simplicity, we shall treat only the affine case---the case of
central, generic arrangements is similar.

Let $\A$ be a generic affine arrangement of $n$ hyperplanes in $\C^{\ell}$, 
with $n>\ell>1$.   The cohomology algebra $A=H^*(X,\Q)$ is the quotient
of the exterior algebra on generators $e_1,\dots ,e_n$ in degree $1$ by the
ideal $I$ generated by all monomials of the form $e_{i_1}\cdots e_{i_{\ell+1}}$.
Since $I$ is generated in degree $\ell+1>2$, the algebra $A$ is not quadratic,
and hence, not Koszul.  Consequently, if $Y=M(\A^{k+1})$, then,
by Theorem~\ref{thm:RF1},
\begin{equation}
\label{eq:nonRF}
E_*(Y) \not\cong L_*(X)[k]\, ,
\end{equation}
even as graded vector spaces.

To understand in a concrete fashion the reason
for this failure, we first
establish an analogue of a well-known theorem
of Hattori~\cite{Ha}, in the setting of redundant generic arrangements.

Let $\A=\{H_1,\dots ,H_n\}$ be an arbitrary hyperplane arrangement 
in $\C^{\ell}$, and let $\A^k=\{H_1^{\times k},\dots ,H_n^{\times k}\}$ 
be the associated subspace arrangement in $\C^{k\ell}$.  Suppose that 
$H_i$ is defined by $\{ \a_i=0\}$. Then recall from \cite{CCX} that 
$H_i^{\times k}$ is defined by $\{\a_i^{\times k}=0\}$, where 
$\a_i^{\times k}\co \C^{k\ell}\to \C^k$ is the map which sends 
$z$ to $(\a_i(z_{11},\dots ,z_{1\ell}),\dots,  
\a_i(z_{k1},\dots ,z_{k\ell}))$.
For $z\in \C^{k\ell}\setminus H_i^{\times k}$, set
$\tilde\a_i^{\times k}(z)=\a_i^{\times k}(z)/
\norm{\a_i^{\times k}(z)}\in S^{2k-1}$.
Now define a map
\begin{align*}
f^{k}_{\A} : M(\A^k)\subset \C^{k\ell} &\to 
\big(S^{2k-1}\big)^{\times n}\subset \C^{kn}\\
z &\mapsto (\tilde\a_1^{\times k}(z),\dots ,\tilde\a_n^{\times k}(z)) .
\end{align*}

In \cite{Ha}, Hattori showed that the complement of a generic affine
arrangement of $n$ hyperplanes in $\C^{\ell}$ ($n\ge \ell$)
is homotopy equivalent to the
$\ell$--skeleton of the $n$--torus $(S^1)^{\times n}$.
The next Lemma extends this result to the corresponding
redundant subspace arrangements.

\begin{lemma}
\label{lem:fatwedge}
Suppose $\A$ is a generic affine arrangement of $n$ hyperplanes
in $\C^{\ell}$, with $n \geq \ell$. Let 
$(S^{2k+1})^{\times n}_{\ell}$ be the $(2k+1)\ell$--skeleton of 
the standard CW--decomposition of $(S^{2k+1})^{\times n}$.
Then, for each $k\ge 1$, the map
$f^{k+1}_{\A}\co M(\A^{k+1})\to \big(S^{2k+1}\big)^{\times n}$
factors  through a homotopy equivalence
\[
f^{k+1}_{\A}\co M(\A^{k+1})\xrightarrow{\:\simeq\:} 
(S^{2k+1})^{\times n}_{\ell}.
\]
\end{lemma}

\begin{proof}
Set $X=M(\A)$ and $Y=M(\A^{k+1})$. From \cite{Ha},
we know that $X\simeq  (S^{1})^{\times n}_{\ell}$.
From \cite[Corollary 2.3]{CCX}, we find that
$Y$ has the same integral homology as
the skeleton $T=(S^{2k+1})^{\times n}_{\ell}$.

Now, since $Y$ is $1$--connected,
it admits a minimal cell decomposition, see \cite{An}.  In particular,
$Y$ is a CW--complex of the same dimension as $T$.
Hence, by cellular approximation, the map $f^{k+1}_{\A}$
factors through $T$. Clearly, the map
$f^{k+1}_{\A}\co Y\to T$
induces an isomorphism in homology.
The result follows from Whitehead's theorem.
\end{proof}

As is well-known, fat wedges of spheres are not coformal
(see eg~\cite{Ta}). Thus, the
reason formula~\eqref{eq:RF} fails for an affine,
generic arrangement $\A$ (with $n>\ell>1$)
is that each of the redundant subspace arrangements,
$\A^{k+1}$ ($k\ge 1$), has non-coformal complement.
If $\ell=n-1$, the obstructions to coformality are
precisely the higher-order Whitehead products, which,
in this case, account for the deviation from equality
in the Rescaling Formula.

\begin{example}
\label{ex:threelines}

Let $\A$ be an affine, generic arrangement of $n$ hyperplanes in 
$\C^{n-1}$, $n>2$. In this case, $X\simeq (S^{1})^{\times n}_{n-1}$, 
and $Y\simeq (S^{2k+1})^{\times n}_{n-1}$. Obviously,
$L_*(X)[k]=\BL^{\ab}(x_1,\dots,x_n)$, with $\deg x_i=2k$.
Using a result of Porter~\cite{Po} on homotopy groups of fat wedges,
we find:
\begin{equation}
\label{eq:threelines}
E_*(Y)=\BL^{\ab}(x_1,\dots,x_n) \amalg \BL(w),
\end{equation}
where $w$ corresponds to the $n$--fold Whitehead product
of $x_1,\dots ,x_n$,  and $\amalg$ stands for the free product
of Lie algebras.  The deviation from equality in the
Rescaling Formula first appears in $\deg w=(2k+1)n-2$, where we
have $\dim E_{(2k+1)n-2}(Y)=1$, but  $\dim L_{(2k+1)n-2}(X)[k]=0$.
Note also that
\[
P_{\Omega Y}(t)=\frac{1}{(1-t^{2k})^{n}- t^{(2k+1)n-2}}\, ,
\]
whereas the homotopy LCS formula would predict
$$P_{\Omega Y}(t)=\frac{1}{(1-t^{2k})^{n}+(-1)^{n+1} t^{2kn}}.$$
\end{example}

\section{Rescaling spherical links}
\label{sec:links}

In this section, we show how to perform geometrically the rescaling operation
for links of circles in $S^3$, and we prove Theorem~\ref{thm:RFlk}
(and its corollary) from the Introduction.

\subsection{Joins of links}
\label{subsec:joinlinks}

A ({\em spherical}) {\em $p$--link} in $S^{2p+1}$ ($p\ge 1$)
is an ordered collection, $K=(K_1,\dots ,$ $ K_n)$, of
pairwise disjoint, smoothly embedded $p$--spheres in the $(2p+1)$--sphere.
In other words, a $p$--link is the image of a smooth embedding 
$K\co \coprod^n S^p \to S^{2p+1}$.
A $1$--link is just a (classical) link.

Recall that the join of two spaces, $X$ and $Y$, is the union,
$X*Y=CX\times Y\cup X\times CY$, glued along $X\times Y$,
where $CX$ denotes the cone on $X$.
For example, $S^p*S^{p'}=S^{p+p'+1}$.
The join operation is functorial:  to maps $f\co X\to X'$
and $g\co Y\to Y'$ there corresponds a map of joins,
$f*g\co X*Y\to X'*Y'$.

As noted by Koschorke and Rolfsen \cite{KR},
the join construction has an analogue for links.

\begin{definition}\cite{KR}\qua
\label{def:join}
Let  $K=(K_1,\dots ,K_n)$ be a $p$--link, and
$K'=(K'_1,\dots ,$ $K'_{n})$ be a $p'$--link
(with the same number of components).
Their {\em join}\/ is the $(p+p'+1)$--link
\[
K\circledast K'=(K_1*K'_1,\dots, K_n*K'_n).
\]
\end{definition}
More precisely, if the two links are the images of embeddings
$K\co \coprod^n S^p \to S^{2p+1}$ and
$K'\co \coprod^n S^{p'} \to S^{2p'+1}$, then their
join is the image of the embedding
$K\circledast K'=\coprod^n K_i*K_i':
\coprod^n S^p*S^{p'}\to S^{2p+1}*S^{2p'+1}$.

\begin{example}
\label{ex:hopf}
Let $\H_n=\widehat{\Delta_n^2}$ be the $n$--component {\em Hopf link},
defined as the Artin closure of the full-twist braid on $n$ strings
(see Figure~\ref{fig:hopfborro} for $n=4$).
As shown in \cite{O}, the Hopf link is the singularity link of
$\A_n=\{z_1^n-z_2^n=0\}$, the arrangement of $n$ complex
lines through the origin of $\C^2$.  It follows from \cite{Ne} 
that $\H_n \circledast \H_n$ is the singularity link of
the subspace arrangement $\A_n^2=\{z_1^n-z_2^n=z_3^n-z_4^n=0\}$.
\end{example}

\subsection{Linking numbers}
\label{subsec:lk}
Let $K=(K_1,\dots ,K_n)$ be a $p$--link in $S^{2p+1}$.
Fix orientations on the ambient sphere
$S^{2p+1}$,  and on each component of $K$.
It is known that all components of $K$ have trivial normal bundle
in $S^{2p+1}$; see Kervaire \cite[Theorem 8.2]{Ke}.
Denote by $T_i\cong S^p \times S^p$ the boundary of a tube,
$K_i\times D^{p+1}$, around $K_i$.
Let $X=S^{2p+1}\setminus \bigcup_{i=1}^{n} K_i$ be the complement
of the link, and let
$\overline{X}=S^{2p+1}\setminus \bigcup_{i=1}^{n} K_i\times
\operatorname{int} D^{p+1}$
be the exterior. We will need some classical definitions and results,
related to the description of the cohomology ring of the complement
in terms of linking numbers.

By standard duality arguments, we know that both $\widetilde{H}_*(X, \Z)$ and
$\widetilde{H}^*(X, \Z)$ are finitely-generated, free abelian groups,
concentrated in two degrees, $p$ and $2p$. Canonical bases may be constructed
as follows. The group $H_p(X, \Z)$ is freely generated by the orientation
classes of the meridional spheres. Denote by $\{ a_1,\dots,a_n \}$ the
Kronecker-dual basis of $H^p(X, \Z)$. For each $i \neq j$, denote by
$b_{ij}\in H^{2p}(X, \Z)$ the Lefschetz dual of an (embedded) arc which
connects, in the exterior of $K$, the tubes $T_i$ and $T_j$.
It is easy to check that
$\{b_{ij}\}$ is a generating set for $H^{2p}(X, \Z)$, with relations
$b_{ij}+b_{ji}=0$ and $b_{ij}+b_{jk}+b_{ki}=0$.

The description of the multiplication in $H^*(X, \Z)$ involves {\em linking
numbers}. For each $i\neq j$, denote by $l_{ij} :=\lk(K_i, K_j)\in \Z$ the
orientation class of $K_i$ in $H_p(S^{2p+1}\setminus K_j, \Z)$, with respect 
to the meridional basis element. With the above notation, one has
the following basic formula:
\begin{equation}
\label{eq:aicupaj}
a_i \cup a_j = (-1)^{p+1}l_{ij}b_{ij}\, , \quad \text{for all} \,\, i\neq j\,.
\end{equation}

To derive \eqref{eq:aicupaj}, it is enough to evaluate $a_i \cup a_j$ 
on $[T_r]$, for each $r$, and to prove that the result equals either
$(-1)^{p+1}l_{ij}$ (if $r=j$), or $(-1)^{p}l_{ij}$ (if $r=i$),
or $0$ (otherwise).  This, in turn, follows from standard cup-product 
computations in the torus $T_r=S^p\times S^p$, given the fact that 
$l_{ij}$ is, by construction, the coefficient of the $i$-th longitude 
on the $j$-th meridian.

\subsection{Linking numbers and joins}
\label{subsec:lkjoin}
Viewing the link $K$ as an embedding in $\R^{2p+1}$,
the linking number $l_{ij}$ equals the degree of the mapping
\[
S^p\times S^p \to S^{2p}, \quad (x,y)\mapsto
(K_i(x)-K_j(y))/ \norm{K_i(x)-K_j(y)},
\]
see \cite[Theorem 5.D.2]{Ro} for $p=1$, and \cite[Lemma 5.1]{Ke} 
for $p>1$.

In view of this remark, we may apply Proposition~4.9 from Koschorke
and Rolfsen~\cite{KR}  to conclude that the linking numbers multiply
under the join operation.  More precisely, if $K=(K_1,\dots ,K_n)$
is a $p$--link and $K'=(K'_1,\dots,K'_n)$ is a $p'$--link, then
\begin{equation}
\label{eq:starlink}
\lk(K_i*K_i',K_j*K_j')=\lk(K_i,K_j)\lk(K_i',K_j').
\end{equation}

\subsection{Rescaling links}
\label{subsec:RFlink}
We are now ready to define the rescaling operation for classical links,
and show that the complement of the rescaled link is indeed the
rescaling of the link complement.

\begin{definition}
\label{def:linkrescale}
Let $K=(K_1,\dots ,K_n)$ be a link in $S^3$.
For each $k\ge 1$, the {\em $k$--rescaling} of $K$
is the following $(2k+1)$--link in $S^{4k+3}$:
\begin{equation}
\label{eq:starhopf}
K^{\circledast k}=K\circledast \underbrace{\H_n \circledast\,
\cdots \circledast \H_n}_{k\ \text{times}}.
\end{equation}
\end{definition}

For example, $\H_n^{\circledast k}$ is the singularity link of
the redundant subspace arrangement $\A_n^{k+1}$; thus, its
complement is homotopy-equivalent to  $M(\A_n^{k+1})$.

\begin{prop}
\label{prop:reslice}
Let $X$ be the complement of a classical link, $K$.
For each $k\ge 1$, the complement, $Y$, of the link $K^{\circledast k}$
is an (integral) $k$--rescaling of $X$.  Moreover, $Y$ is
the unique (up to $\Q$--equivalence) $k$--rescaling of $X$;
in particular, $Y$ is a formal space.
\end{prop}

\begin{proof}
Since $K^{\circledast k}$ has codimension $2(k+1)>2$,
the complement $Y$ is simply-connected.  
Let $a_i\in H^1(X,\Z)$ and $\tilde{a}_i\in H^{2k+1}(Y,\Z)$ be the 
standard basis elements (Kronecker dual to the meridional spheres). 
Since  $a_i$ has odd degree, and $H^{*}(X,\Z)$ is torsion-free, 
we infer that $a_i^2=0$, for all $i$; similarly, $\tilde{a}_i^2=0$. 
Since all the linking numbers of the Hopf link are equal to $1$,
and since both $K$ and $K^{\circledast k}$ are odd-dimensional links,
we see from equations \eqref{eq:aicupaj}, \eqref{eq:starlink},
and \eqref{eq:starhopf} that
\[
H^*(Y,\Z)=H^*(X,\Z)[k], \quad\text{as graded rings}.
\]
Therefore, $Y$ is a $k$--rescaling of $X$, even integrally.

Since $H^{>2}(X,\Z)=0$, the uniqueness up to rational
homotopy (and thereby, the formality) of
the rescaling follows from Proposition~\ref{sss1}.
\end{proof}

\subsection{Proofs of Theorem~\ref{thm:RFlk} and Corollary~\ref{cor:lkresc}}
\label{subsec:pfE}
Results from \cite{MP} (Theorem B'(i) and the Corollary to Proposition 6.3)
show that $H^*(X, \Q)$ is a Koszul algebra if and only if the linking graph
$\G_K$ is connected. Theorem~\ref{thm:RFlk} follows then from Theorem
\ref{thm:RF1}, via Remark~\ref{rem:asymp} and Proposition~\ref{prop:reslice}.

Assume now that $\G_K$ is connected. The coformality of $Y$ is 
ensured by Proposition~\ref{prop:coformal}. The semidirect product 
structure of the Lie algebra $E_*(Y)$ follows, via the Rescaling Formula, 
from the corresponding semidirect product decomposition of $L_*(X)=\HH_*$ 
(see Lemma~\ref{lem:r1}), proved in \cite[Theorem 4.2(i) and Lemma 4.1]{BP}.  
The other statements in Corollary~\ref{cor:lkresc} are direct consequences 
of Theorem~\ref{thm:HLCS} and the discussion in \S\ref{subsec:ratloop}.

\section{Formality and the Campbell-Hausdorff invariants}
\label{sec:CH}

We start by reviewing the Campbell--Hausdorff invariants of links,
introduced in \cite{P97}. We then explore the relationship between
these invariants and the formality of the link complement.
We conclude with an example showing that the Malcev Formula~\eqref{eq:Mfinal}
from Theorem~\ref{thm:MF} may fail for a non-formal space $X$, even if
$H^*(X, \Q)$ is a Koszul algebra.

\subsection{CH-invariants}
\label{subsec:chinv}
Let $\BL_n=\BL^*(x_1,\dots,x_n)$ be the free Lie algebra on 
$n$ generators $x_i$ in degree~$1$, and let $\wL_n$ be its completion
with respect to bracket length, constructed as in Example~\ref{ex:mformal}.
Let $V$ be an $n$--dimensional $\Q$--vector space,
with basis $v_1,\dots ,v_n$. The {\em universal moduli space}
for links of $n$ components is the $\Q$--vector space
\[
\Der :=\Hom_{\Q}(V, F_2\wL_n) \, ,
\]
with elements, $\partial \in \Der$, written in formal series notation,
\[
\partial = \sum_{i\ge 1} \partial_i\, , \quad \text{where} \quad
\partial_i \in \Hom_{\Q}(V,\BL_n^{i+1}) \, .
\]
The vector space $\Der$ comes endowed with a filtration,
defined by
\[
F_r\Der=\{ \partial \in \Der \mid \partial_i= 0, \text{ for $i<r$}\}\, .
\]
There is then a certain pro-unipotent $\Q$--group, $U_n$,
which acts linearly on $\Der$, so as to preserve the filtration
$\{F_r\Der\}_{r\ge 1}$.

Let $\F_n$ be the free group on $x_1,\dots,x_n$.
Let $\exp(\wL_n)$ be the exponential group corresponding
to the Malcev Lie algebra $\wL_n$.  The Campbell--Hausdorff 
representation, $\rho\co \F_n \to \exp(\wL_n)$, is defined by 
$\rho(x_i)=x_i$.

Now let $K=\{K_1,\dots ,K_n\}$ be a link in
$S^3$, with a fixed ordering and orientation of the components.
Attached to $K$ there is a homomorphism $\partial=\partial_{K,\{l_i\}} \in \Der$,
constructed by Campbell--Hausdorff expansion from the
longitudes $l_i\in \F_n$ of the link $K$:
\begin{equation}
\label{eq:initial}
\partial (v_i):=\left[x_i,\rho(l_i)\right]=\sum_{j=1}^n
\left[ x_i,l_{ij}x_j\right] +\cdots
\end{equation}
Note that the degree $1$ part, $\partial_1$, is determined solely
by the linking numbers, $l_{ij}$, of $K$.

\begin{definition}\cite{P97}\qua
\label{def:ch}
The {\em Campbell--Hausdorff invariant of order $r$} of the link $K$ is
the $U_n$--orbit of $\partial$ in the quotient of
$\Der$ by $F_r\Der$:
\[
p^r(K)=\overline{\partial} \in U_n\backslash \Der/F_r\Der \, .
\]
\end{definition}
It turns out that $p^r(K)$ depends only on $K$, and not on the
choice of longitudes; in fact, $p^r(K)$ depends only on
the concordance class of $K$; see \cite{P97}.

\subsection{Relation with formality}
\label{subsec:CHformal}
The next results make precise the connection between the CH-invariants
of a link and the formality of the link complement.

\begin{lemma}
\label{lem:obstr}
Let $K$ be a link with connected linking graph. Choose any system of longitudes,
$\{l_i\}$, and associate to it an element $\partial=\partial_{K,\{l_i\}} \in \Der$,
as in \eqref{eq:initial}. If the complement, $X$, of $K$ is a formal space, then
$\partial$ and $\partial_1$ lie in the same $U_n$--orbit.
\end{lemma}

\begin{proof}
This follows from the proof of \cite[Corollary 3.6]{P01}. The
connectedness of $\G_K$ comes into play, via the Koszul property
of the algebra $A^* =H^*(X, \Q)$, to guarantee that the comultiplication map,
$\nabla \co A_2 \to A_1\wedge A_1$, is injective. (See \cite[Remark 3.4]{P01}
for the connection between $\partial_1$ and $\nabla$.) The point is that the
injectivity of $\nabla$, together with the formality of $X$, are actually
sufficient for the proof of Corollary 3.6 from \cite{P01}, which gives
the desired conclusion.
\end{proof}

\begin{cor}
\label{cor:obstr}
Let  $K$ and $K_0$ be two links with the same connected
weighted linking graph.  If both link complements are formal,
then $p^r(K)=p^r(K_0)$, for all $r\ge 1$.
\end{cor}

\begin{proof}
This follows from the Lemma above, given the
definition~\ref{def:ch} of the CH-invariants, and the fact that
$\partial_1$  is determined by the linking numbers. 
\end{proof}

\subsection{A non-formal example}
\label{subsec:nonMF}
Using the Campbell--Hausdorff invariants, we exhibit examples
of link complements with Koszul cohomology algebra, for which the
Malcev Formula does not hold.

\begin{figure}
\setlength{\unitlength}{1cm}
\label{fig:hopf}
\begin{minipage}[t]{0.8\textwidth}
\begin{picture}(6,1.8)(-4,-0.5)
\xy
0;/r1pc/:
,{\htwist\htwist}
,+(0,1)
,{\htwist}
,+(0,-1)
,{\htwist\htwist}
,+(0,1)
,{\htwist}
,+(0,1)
,{\htwist}
,+(0,-1)
,{\htwist}
,+(0,-1)
,{\htwist\htwist}
,+(0,1)
,{\htwist}
,+(0,1)
,{\htwist}
,{\htwist~{(3,-1)}{(4,-1)}{(3,-1)}{(4,-1)}}
,{\htwist~{(6,-1)}{(9,-1)}{(6,-1)}{(9,-1)}}
,{\htwist~{(11,-1)}{(13,-1)}{(11,-1)}{(13,-1)}}
,{\htwist~{(7,0)}{(8,0)}{(7,0)}{(8,0)}}
,{\htwist~{(12,0)}{(13,0)}{(12,0)}{(13,0)}}
,{\htwist~{(1,1)}{(3,1)}{(1,1)}{(3,1)}}
,{\htwist~{(4,1)}{(6,1)}{(4,1)}{(6,1)}}
,{\htwist~{(9,1)}{(11,1)}{(9,1)}{(11,1)}}
,{\htwist~{(1,2)}{(7,2)}{(1,2)}{(7,2)}}
,{\htwist~{(8,2)}{(12,2)}{(8,2)}{(12,2)}}
\endxy
\end{picture}
\end{minipage}

\label{fig:hopfb}
\begin{minipage}[t]{0.8\textwidth}
\begin{picture}(6,2.5)(-2.4,-0.8)
\xy
0;/r1pc/:
,{\htwist\htwist}
,+(0,1)
,{\htwist}
,+(0,-1)
,{\htwist\htwist}
,+(0,1)
,{\htwist}
,+(0,1)
,{\htwist}
,+(0,-1)
,{\htwist}
,+(0,-1)
,{\htwist\htwist}
,+(0,1)
,{\htwist}
,+(0,1)
,{\htwist}
,+(0,-2)
,{\htwist}
,+(0,1)
,{\htwist,\htwist}
,+(0,-1)
,{\htwistneg,\htwistneg}
,+(0,1)
,{\htwistneg,\htwistneg}
,+(0,-1)
,{\htwist}
,{\htwist~{(3,-1)}{(4,-1)}{(3,-1)}{(4,-1)}}
,{\htwist~{(6,-1)}{(9,-1)}{(6,-1)}{(9,-1)}}
,{\htwist~{(11,-1)}{(13,-1)}{(11,-1)}{(13,-1)}}
,{\htwist~{(7,0)}{(8,0)}{(7,0)}{(8,0)}}
,{\htwist~{(12,0)}{(13,0)}{(12,0)}{(13,0)}}
,{\htwist~{(1,1)}{(3,1)}{(1,1)}{(3,1)}}
,{\htwist~{(4,1)}{(6,1)}{(4,1)}{(6,1)}}
,{\htwist~{(9,1)}{(11,1)}{(9,1)}{(11,1)}}
,{\htwist~{(1,2)}{(7,2)}{(1,2)}{(7,2)}}
,{\htwist~{(8,2)}{(12,2)}{(8,2)}{(12,2)}}
,{\htwist~{(14,-1)}{(16,-1)}{(14,-1)}{(16,-1)}}
,{\htwist~{(18,-1)}{(20,-1)}{(18,-1)}{(20,-1)}}
,{\htwist~{(13,1)}{(14,1)}{(13,1)}{(14,1)}}
,{\htwist~{(16,1)}{(18,1)}{(16,1)}{(18,1)}}
,{\htwist~{(20,1)}{(21,1)}{(20,1)}{(21,1)}}
,{\htwist~{(13,2)}{(21,2)}{(13,2)}{(21,2)}}
\endxy
\end{picture}
\end{minipage}
\caption{Braids $\Delta_4^2$ and
$\Delta_4^2\gamma$}
\label{fig:hopfborro}
\end{figure}
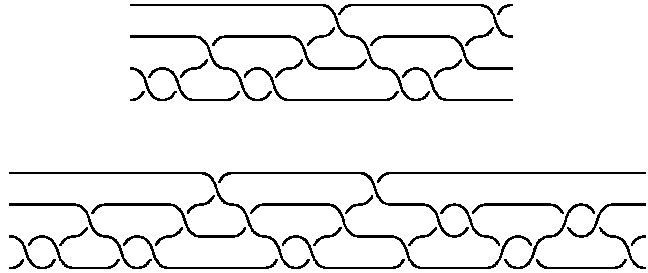

\begin{example}
\label{ex:twisthopf}

Let $K_0=\widehat{\Delta_n^2}$ be the $n$--component Hopf link.
Since $K_0$ is the singularity link of the complex line arrangement $\A_n$,
its complement, $X_0$, is formal.

Assume $n>3$. Pick any pure braid $\gamma$ in
$\Gamma_{r-1} P_{n-1}\setminus \Gamma_{r} P_{n-1}$,
with $r\ge 3$, and let $K=\widehat{\Delta_n^2\gamma}$ be the
link obtained by closing up the full-twist braid followed by $\gamma$.
(See Figure~\ref{fig:hopfborro} for a simple example,
with $n=4$ and $r=3$.)  Let $X$ be the complement of $K$.

The links $K$ and $K_0$ are not that easy to distinguish.
Indeed, they both have:
\begin{itemize}
\item
The same weighted linking graph: the complete graph on $n$ vertices,
with all linking numbers equal to~$1$. (This uses the fact that $r\ge 3$.)
Thus, $H^*(X,\Z)\cong H^*(X_0,\Z)$, as graded rings.

\item
The same Milnor $\overline{\mu}$--invariants, of arbitrary length
(they all vanish, except for the linking numbers).

\item
The same Vassiliev invariants, up to order $r-2$; see \cite[Theorem 1]{St}.

\item
The same  associated graded {\em integral} Lie algebras; 
see \cite[Corollary 6.2]{MP}.
\end{itemize}

On the other hand, as shown in \cite[Proposition~6.1]{P01}, the
Campbell--Hausdorff invariants of order $r$ distinguish the two links:
\[
p^{r}(K)\ne p^{r}(K_0).
\]
Hence, by Corollary~\ref{cor:obstr}, the space $X$ is non-formal.

Let $Y$ be a homological $k$--rescaling of $X$ (for example,
the complement of $K^{\circledast k}$).  
We know that $H^*(X, \Q)$ is a Koszul algebra,
since $\G_K$ is connected.  Hence, by Theorem~\ref{thm:RF2},
the Rescaling Formula holds:  $E_*(Y)=L_*(X)[k]$.
On the other hand, we also know that $X$ is not formal.
Hence, by  Corollary~\ref{cor:Mform},
the Malcev Formula fails in this case:
\[
\pi_1(X)\otimes \Q\not\cong\comp{[\Omega S^{2k+1}, \Omega Y]} \quad
\text{(as filtered groups)}.
\]
\end{example}


\end{document}